\documentclass{article}
\usepackage{graphics}
\usepackage{graphicx}
\usepackage{subcaption}
\usepackage{bm}
\usepackage{xcolor}
\usepackage{amsmath}
\usepackage{longtable}
\usepackage{multirow}
\usepackage{booktabs}
\usepackage{amssymb}
\usepackage{amsmath}
\usepackage{listings}  
\usepackage{xcolor}    
\usepackage{caption}   

\begin{document}
\parindent=22pt
\normalsize \baselineskip=18pt
\date{}
\renewcommand{\baselinestretch}{1.2}
\renewcommand{\arraystretch}{1.5}
\renewcommand{\thefootnote}{\fnsymbol{footnote}}
\begin{center}
{\Large\bf Statistical inference of partially linear time-varying coefficients spatial autoregressive panel data model}
\vskip 0.5cm
{\large  Lingling Tian$^1$ \quad  Chuanhua Wei$^2$ \quad Mixia Wu
$^1$ \footnote{ Corresponding Author:
  Mixia Wu, School of Mathematics, Statistics and Mechanics, Beijing
 University of Technology, Beijing  100124, P.R. China;
 E-mail:wumixia@bjut.edu.cn}}
 \vskip 0.4cm
 1.School of Mathematics, Statistics and Mechanics,
 Beijing University of Technology, Beijing 100124, P.R. China \\
 2.School of Science, Minzu
 University of China, Beijing 100081, P.R. China \\
  \vskip 1cm
{\large\bf Abstract}
\end{center}
\large

This paper  investigates  a partially linear spatial autoregressive panel data model that incorporates fixed effects, constant and time-varying regression coefficients, and a time-varying spatial lag coefficient. A two-stage least squares estimation method based on profile local linear dummy variables (2SLS-PLLDV) is proposed to estimate both constant and time-varying coefficients without the need for first differencing. The asymptotic properties of the estimator are derived under certain conditions. Furthermore, a residual-based goodness-of-fit test is constructed for the model, and a residual-based bootstrap method is used to obtain p-values. Simulation studies show the good performance of the proposed method in various scenarios.  The Chinese provincial carbon emission data set is analyzed for illustration.
\\


\section{Introduction}
Spatial panel data are prevalent in fields such as geography, regional economics, environmental science, and public health. In recent years, the spatial autoregressive (SAR) panel data model has gained a great deal of attention due to its ability to account for spatial dependence and temporal dynamics simultaneously. Various spatial autoregressive panel data models have been proposed, including the classical linear SAR model (Lee 2010), semiparametric and nonparametric SAR models (see, Ai and Zhang 2015; Tian RQ {\it et al.} 2024; Chen {\it et al.} 2019; Liang {\it et al.} 2022). The estimation methods for SAR panel data models have undergone rapid developments. For instance, Lee (2010) proposed the quasi-maximum likelihood (QML) estimation methods for linear SAR panel data model. For the partially linear SAR panel data model, Xu and Chen (2017) built the profile quasi-maximum likelihood (PQML) estimation, while Ai and Zhang (2015) developed a sieve two stage least squares (S2SLS) method. Liang {\it et al.} (2022) introduced a local linear quasi-maximum likelihood (LLQML) estimation method for a semiparamertic SAR model. Notably, the coefficient of the spatial lag term in these corresponding models  remains constant.
   
To capture more complex autoregressive relationships, Sun and Malikov (2018) studied a functional-coefficient SAR model, which allows the coefficient of the spatial lag term to vary with covariates. Extending this concept further, Chang {\it et al.} (2024) considered a trending time-varying coefficient SAR panel data model in which  all unknown coefficients, including  the spatial lag term coefficient, are allowed to vary over time. The model can be given by
\begin{equation}   \label{1}
y_{i t}=\rho(\tau_t) \sum_{j\neq i} w_{i j} y_{j t}+{\bf
	x}_{it}^{\top}{\boldsymbol\beta}(\tau_{t})+\alpha_{i}+\varepsilon_{i t}, ~~i=1,2, \cdots, N,~ t=1,2, \cdots, T, \end{equation}
where $\rho(\tau_t)$ is the time-varying coefficient of spatial lag term $\sum_{i\neq j}w_{i j} y_{j t}$, constrained by the condition $|\rho(\tau_t)|<1$.  The $w_{i j}$s constitute the elements of the pre-determined spatial adjacency matrix ${\bf W}_{N}$, which describes the spatial distance between location $i$ and $j$,     with $w_{ii}=0$  explicitly stating that there is no self-connection.  The matrix ${\bf W}_{N}$ is commonly  row-standardized. The observations $y_{it}$ and ${\bf x}_{it}=(x_{it1},x_{it2},\cdots,x_{itp})^{\top}$ correspond to the response variable   $ Y$ and $p$-dimensional explanatory variables ${\bm X}=(X_1,X_2,\cdots,X_p)^\top$ respectively, at the $i$-th location and the $t$-th time.  Here, $N$ and $T$ denote the total number of spatial units and the length of the time series, respectively. The vector ${\boldsymbol\beta}(\tau_t)=\left(\beta_{1}(\tau_{t}),\beta_{2}(\tau_{t}),\cdots,\beta_{p}(\tau_{t})\right)^{\top}$ 
contains the $p$ unknown coefficients associated with the explanatory variables ${\bf x}_{it}$, where $\tau_{t}= {t}/{T} \in(0,1]$ represents a time-normalized index. Conventionally, $x_{it1}=1$ is assumed to enable the model to incorporate a time-varying intercept (it is referred to as the "time trend" in Chang {\it et al.} (2024)). For the purpose of identifying ${\bm \beta}(\tau_t)$, the individual fixed effects $\alpha_{i}$s are constrained  by the condition $\sum_{i=1}^N\alpha_{i}=0$. The error terms $\varepsilon_{i t}$s are assumed to be independent and identically distributed with zero mean and a finite variance $\sigma_{\bm \varepsilon}^2 (>0)$.
 

 \par

 
It is worth noting that as $p$ increases, despite its inherent flexibility,  Model (\ref{1})  encounters the curse of dimensionality in estimating the $p$ unknown coefficient functions. 
If some coefficients of  $\{\beta_{1}(\tau_{t}),\beta_{2}(\tau_{t}),\cdots,\beta_{p}(\tau_{t}\}$ are constant over time,  
 for simplicity and without loss of generality, 
 we assume that ${\beta}_k(\tau_{t}) \equiv {\beta}_k ,~ k=q+1, \cdots,p,$ 
then we obtain a partially linear time-varying coefficients SAR panel data models with fixed effects: 
\begin{equation}    \label{2}
    y_{i t}=\rho(\tau_t) \sum_{j\neq i} w_{i j} y_{j t}+{\bf
	x}_{v,it}^{\top}{\boldsymbol\beta}_v(\tau_{t})+{\bf
	x}_{c,it}^{\top}{\boldsymbol\beta}_c+\alpha_{i}+\varepsilon_{i t}, i=1,2, \cdots, N; t=1,2, \cdots, T, 
\end{equation}
where  
${\bf x}_{v,it}=(x_{it1},x_{it2},\cdots,x_{itq})^{\top}$ 
comprises the first $q$ observations,  
$\boldsymbol\beta_v(\tau_t)=(\beta_1(\tau_t),$ $\beta_2(\tau_t),\cdots,\beta_q(\tau_t))^{\top}$ is the corresponding  vector of the time-varying coefficients, and ${\boldsymbol\beta}_c=\left(\beta_{q+1},\beta_{q+2},\cdots,\beta_{p}\right)^{\top}$ represents the vector of  the time-invariant coefficients of the remaining $p-q$ observations ${\bf x}_{c,it}=(x_{it,q+1},x_{it,q+2},\cdots,x_{itp})^{\top}  $. It is evident that when $q=p$, Model (\ref{2}) simplifies to Model (\ref{1}), indicating that Model (\ref{2}) offers greater flexibility than Model (\ref{1}). 

To the best of our knowledge, there is limited literature on Model (\ref{2}). This paper primarily aims to estimate Model (\ref{2}) by proposing a two-stage least squares estimation method based on profile local linear dummy variables, referred to as the 2SLS-PLLDV method. 
   The 2SLS method relies on selecting instrumental variables in the first stage to obtain consistent estimators of endogenous explanatory variables. It is a critical issue for how to choosing effective instrumental variables. Kelejian and Prucha (1998), who applied the 2SLS method to SAR models, recommended constructing instrumental variables using linearly independent columns such as  $({\bf X}, {\bf WX}, {\bf W}^2{\bf X},\cdots)$, where ${\bf X}$  represents the design matrix and ${\bf W}$ denotes the spatial weighted matrix, respectively. Following this,  we  develop  
   a goodness-of-fit test for Model (2)    against the alternative hypothesis represented by Model (1), specifically to test the following hypothesis:
  \begin{equation}\label{3}
{H_{0}}: {\beta}_k(\tau_{t}) \equiv {\beta}_k ,~ k=q+1, \cdots,p  \quad {\rm{VS}} \quad   H_{1}: {\beta}_k(\tau_{t}) \neq {\beta}_k,~ k=1, \cdots,p. 
\end{equation}

\par
The structure of this paper is outlined as follows. In Section 2, we propose the  2SLS-PLLDV  estimation for Model (\ref{2}). 
In Section 3, we establish the asymptotic properties of the proposed estimators. In Section 4, we  constructs a test statistic for Hypothesis (\ref{3}), based on the residual sum of squares. The determination of its p-values is facilitated by a residual-based bootstrap procedure. Section 5 presents simulation studies  designed to assess the finite sample performance of our proposed estimation and testing methodology. In Section 6, we demonstrate the practical application of our methods by analyzing provincial carbon emission data from China. Lastly,  Section 7 concludes the paper by summarizing the key findings and contributions.

\section{The 2SLS-PLLDV estimation }
\par 
Denote ${\bm\rho}_{NT}=\mbox{diag}(\rho(\tau_1),\rho(\tau_2),...,\rho(\tau_T))\otimes {\bf I}_N$, ${\bf W}={\bf I}_{T} \otimes {\bf W}_{N}$,  ${\bf D}={\bf 1}_{T} \otimes({\bf -1}_{N-1}, {\bf I}_{N-1})^{\top}$, ${\boldsymbol\alpha}=(\alpha_{2}, \alpha_{3}, \ldots, \alpha_{N})^{\top}$, where   ${\bf W}_{N}=(w_{i j})$, ${\bf I}_n$  represents the identity matrix of order   $n$, and ${\bf 1}_{m}$ is a $m\times 1$ vector whose elements are all one.  Let
$$
{\bm Y}=\left(
\begin{array}{c}
	y_{11}\\
	\vdots\\
	y_{N1}\\
	\vdots\\
	y_{NT}
\end{array}
\right) ,\quad 
{\bf B}({\bf X}_v,{\bm \beta}_v)=\left(
\begin{array}{c}
{\bf x}_{v,11}^\top{\boldsymbol\beta}_v(\tau_1) \\
\vdots\\
{\bf x}_{v,N1}^\top{\boldsymbol\beta}_v(\tau_1) \\
\vdots\\
{\bf x}_{v,NT}^\top{\boldsymbol\beta}_v(\tau_T)
\end{array}
\right),\quad
{\bf X}_c=\left(
\begin{array}{c}
{\bf x}_{c,11}^\top \\
\vdots\\
{\bf x}_{c,N1}^\top \\
\vdots\\
{\bf x}_{c,NT}^\top
\end{array}
\right),\quad
{\boldsymbol\varepsilon}=\left(
\begin{array}{c}
	\varepsilon_{11} \\
	\vdots\\
	\varepsilon_{N1} \\
	\vdots\\
	\varepsilon_{NT}
\end{array}
\right).
$$
Then Model (2),  expressed in matrix form, is   
\begin{equation}\label{4}
{\bm Y}={\bm \rho}_{NT}{\bf W} {\bm Y}+{\bf B}({\bf X}_v,{\bm \beta}_v)+{\bf X}_c\boldsymbol\beta_c+{\bf D}\boldsymbol\alpha+\boldsymbol\varepsilon.
\end{equation}

\subsection{Instrumental estimation of the spatial lag term}

The instrumental variables ${\bf H}=({\bf X},{\bf W}{\bf X}_{(-1)},{\bf W}^2{\bf X}_{(-1)})$ are commonly recommended for estimating the spatial lag term  ${\bf W}\bm Y$, where ${\bf X}$ is the $NT \times p$ design matrix with all elements in the first column being 1, and ${\bf X}_{(-1)}$ is a $NT \times (p-1)$ matrix obtained by removing the first column of $\bf X$. Denote $y_{w,it}=\sum_{j\neq i} w_{i j} y_{j t}$ and  ${\bf W}\bm Y =\bm Y _{\bf w}=(y_{w,11},\cdots ,y_{w,N1},\cdots ,y_{w,NT})^\top. $  Then   the regression of $\bm Y_{\bf w} $ on $\bf H$ can be expressed as follows:
\begin{equation}\label{5}
y_{wit}={\bm h}_{it}^\top {\bm \eta}(\tau_t)+\psi_{i}+e_{it}, 
\end{equation}
where ${\bm h}_{it}=(1,h_{it,2},\cdots,h_{it,3p-2})^\top$ is the row vector of the instrumental variables $\bf H$, ${\bm \eta}(\tau_t)=(\eta_1(\tau_t),\eta_2(\tau_t),\cdots,\eta_{3p-2}(\tau_t))^\top$ is $(3p-2)$-dimensional unknown time-varying coefficients vector, $\psi_i$ is the fixed effects and they satisfy $\sum\limits_{i=1}^N\psi_{i}=0$, $e_{it}$ is the model error. The matrix form of Model (\ref{5})  is given by
$$
{\bm Y}_{\bf w}={\bf B}({\bf H},{\bm \eta})+{\bf D} \boldsymbol\psi+\bm e,
$$
where ${\boldsymbol\psi}=(\psi_{2}, \psi_{3}, \ldots, \psi_{N})^{\top}$. 

Assume that $\boldsymbol\eta(\tau_t)$ has continuous derivatives of up to the second order. Applying Taylor expansion to $\boldsymbol\eta(\tau_{t})$ at $\tau_{0}$, we have
$$
\boldsymbol\eta(\tau_{t})=\boldsymbol\eta(\tau_{0})+\boldsymbol\eta^{\prime}(\tau_{0})({\tau_{t}-\tau_{0}})+\mathrm{O}((\tau_{t}-\tau_{0})^{2}),
$$
where $\boldsymbol\eta^{\prime}(\cdot)$ is the first derivative of $\boldsymbol\eta(\cdot)$. Define the weight matrix as
$$
{\bf W}_h(\tau_{0})={\rm diag}\left[ K\left(\frac{\tau_{1}-\tau_{0}}{h}\right), \ldots, K\left(\frac{\tau_{T}-\tau_{0}}{h}\right)\right] \otimes {\bf I}_{N},
$$
  where ${\tau_0}\in (0,1]$, $K(\cdot)$ is the kernel function and $h$ is the bandwidth. Let $\bm a=\boldsymbol\eta(\tau_{0})$ and $\bm b=h\boldsymbol\eta^{\prime}(\tau_{0})$. Utilizing the local linear approximation, we can deduce the estimators  of $\boldsymbol\eta(\tau_{0})$ and $h\boldsymbol\eta^{\prime}(\tau_{0})$  by solving the following optimization problem:
\begin{equation}\label{6}
\arg \min\limits_{\bm a,\bm  b\in {\bf R}^{3p-2}} \left[{\bm Y}_{\bf w}-{\bf M}_{\bf H}(\tau_{0})(\bm a^{\top},\bm b^{\top})^{\top}-{\bf D} \boldsymbol\psi\right]^{\top}{\bf W}_h(\tau_{0})\left[{\bm Y}_{\bf w}-{\bf M}_{\bf H}(\tau_{0})(\bm a^{\top},\bm b^{\top})^{\top}-{\bf D} \boldsymbol\psi\right], 
\end{equation}
  where
$$
{\bf M}_{\bf H}(\tau_{0})=\left(\begin{array}{cc}
	{\bm h}_{11}^{\top} & \frac{\tau_{1}-\tau_{0}}{h} {\bm h}_{11}^{\top} \\
	\vdots & \vdots \\
	{\bm h}_{N1}^{\top} & \frac{\tau_{1}-\tau_{0}}{h} {\bm h}_{N1}^{\top} \\
	\vdots & \vdots \\
	{\bm h}_{NT}^{\top} & \frac{\tau_{T}-\tau_{0}}{h} {\bm h}_{NT}^{\top}\\
\end{array}\right).
$$
Differentiating Eq.(\ref{6}) with respect to $\bm \psi$ and equating the result to zero,  we obtain
$$
\hat {\bm \psi}(\tau_0)=({\bf D}^\top {\bf W}_h(\tau_0){\bf D})^{-1}{\bf D}^\top {\bf W}_h(\tau_0)[{\bm Y}_{\bf w}-{\bf M}_{\bf H}(\tau_{0})(\bm a^{\top},\bm b^{\top})^{\top}].
$$
Substituting   $\hat {\bm \psi}(\tau_0)$ back into Eq.(\ref{6}), we derive  the concentrated weight least squares: 
\begin{equation}\label{7}
\begin{array}{ll}
	&[{\bm Y}_{\bf w}-{\bf M}_{\bf H}(\tau_{0})(\bm a^{\top},\bm b^{\top})^{\top}-{\bf D}\hat {\bm \psi}(\tau_0)]^{\top}{\bf W}_h(\tau_{0})[{\bm Y}_{\bf w}-{\bf M}_{\bf H}(\tau_{0})\bm a^{\top},\bm b^{\top})^{\top}-{\bf D}\hat {\bm \psi}(\tau_0)]\\
	=&[{\bm Y}_{\bf w}-{\bf M}_{\bf H}(\tau_{0})(\bm a^{\top},\bm b^{\top})^{\top}]^\top{\bf W}_h^{*}(\tau_0)[{\bm Y}_{\bf w}-{\bf M}_{\bf H}(\tau_{0})(\bm a^{\top},\bm b^{\top})^{\top}],
\end{array} 
\end{equation}
where ${\bf W}_h^{*}(\tau_0)= {\mathcal{K}}^\top(\tau_0){\bf W}_h(\tau_{0}){\mathcal{K}}(\tau_0)$  and ${\mathcal{K}}(\tau_0)={\bf I}_{NT}-{\bf D}({\bf D}^\top {\bf W}_h(\tau_0){\bf D})^{-1}{\bf D}^\top {\bf W}_h(\tau_0)$. For any $\tau_0$, we can yield ${\mathcal{K}}(\tau_0){\bf D}{\bm \psi}=0$ through straightforward calculation. Consequently, the fixed effects term   ${\bf D}{\bm \psi}$ is eliminated in Eq.(\ref{7}). Minimizing the Eq.(\ref{7}) with respect to $\bm a$ and $\bm b$, we obtain the estimator of $\boldsymbol\eta(\tau_{0})$  given by 
$$
\hat{\boldsymbol\eta}(\tau_0)={\bm{\Phi}}_{\bf H}(\tau_0){\bm Y}_{\bf w}, 
$$
where ${\bm{\Phi}}_{\bf H}(\tau_0)=\left({\bf I}_{3p-2}, {\bf 0 }_{(3p-2) \times (3p-2)}\right) \left({\bf M}_{\bf H}^{\top}(\tau_0){\bf W}_h^{*}(\tau_0){\bf M}_{\bf H}(\tau_0)\right)^{-1} {\bf M}_{\bf H}^{\top}(\tau_0){\bf W}_h^{*}(\tau_0)$. We refer to $\hat{\boldsymbol\eta}(\tau_0)$ as the local linear dummy variables (LLDV) estimator. For $\tau_0=\tau_1,\cdots,\tau_T$,   the estimator of ${\bf B}({\bf H},\bm\eta)$  can be constructed as follows:
$$
\hat{\bf B}({\bf H},\bm\eta)={\bf S}_{\bf H}(h){\bm Y}_{\bf w},
$$
where
$$
{\bf S}_{\bf H}(h)=\left(\begin{array}{cc}
	({\bf h}_{11}^\top, {\bf 0}_{1 \times (3p-2)})[{\bf M}_{\bf H}^{\top}(\tau_1) {\bf W}_h^{*}(\tau_1) {\bf M}_{\bf H}(\tau_1)]^{-1} {\bf M}_{\bf H}^{\top}(\tau_1) {\bf W}_h^{*}(\tau_1)\\
	\vdots \\
	({\bf h}_{N1}^\top, {\bf 0}_{1 \times (3p-2)})[{\bf M}_{\bf H}^{\top}(\tau_1) {\bf W}_h^{*}(\tau_1) {\bf M}_{\bf H}(\tau_1)]^{-1} {\bf M}_{\bf H}^{\top}(\tau_1) {\bf W}_h^{*}(\tau_1)\\
	\vdots \\
	\left({\bf h}_{NT}^\top, {\bf 0 }_{1 \times (3p-2)}\right)[{\bf M}_{\bf H}^{\top}(\tau_T) {\bf W}_h^{*}(\tau_T) {\bf M}_{\bf H}(\tau_T)]^{-1} {\bf M}_{\bf H}^{\top}(\tau_T) {\bf W}_h^{*}(\tau_T)\\
\end{array}\right).
$$
Then the estimator of $\bm \psi$ is given by 
$$
\hat{\bm \psi}=({\bf D}^\top {\bf D})^{-1}{\bf D}^\top({\bm Y}_{\bf w}-\hat{\bf B}({\bf H},\bm\eta)).
$$
 Subsequently, we derive the estimator of ${\bm Y}_{\bf w}$ as
$$
\hat{{\bm Y}}_{\bf w}=\hat{\bf B}({\bf H},\bm\eta)+{\bf D}\hat{\bm \psi}.
$$
Upon replacing the spatial lag term  ${\bf W}\bm Y$ in Model (4) with the instrumental estimator  $\hat {\bm Y}_{\bf w}$,  we have completed the first stage of the 2SLS-PLLDV  estimation.

\subsection{ PLLDV estimation of regression coefficients }
Substituting $\hat {\bm Y}_{\bf w}$ into the model (4), we obtain
\begin{equation} \label{8}
    {\bm Y}={\bm \rho}_{NT}\hat {\bm Y}_{\bf w}+{\bf B}({\bf X}_v,{\bm \beta}_v)+{\bf X}_c\boldsymbol\beta_c+{\bf D}\boldsymbol\alpha+\boldsymbol\varepsilon.  
\end{equation}
Let ${\bf z}_{v,it}=({\hat y}_{w,it},{\bf x}^\top_{v,it})^\top$, ${\bm \gamma}_v(\tau_t)=(\rho(\tau_t),{\bm \beta}_v^\top(\tau_t))^\top$ and   
$$
{\bf B}({\bf Z}_v,{\bm \gamma}_v)=\left(
\begin{array}{c}
	{\bf z}_{v,11}^\top\boldsymbol\gamma_v(\tau_1) \\
	\vdots\\
	{\bf z}_{v,N1}^\top\boldsymbol\gamma_v(\tau_1) \\
	\vdots\\
	{\bf z}_{v,NT}^\top\boldsymbol\gamma_v(\tau_T)
\end{array}
\right).
$$
Then Model  (\ref{8}) can be rewritten as
\begin{equation} \label{9}
  {\bm Y}={\bf B}({\bf Z}_v,{\bm \gamma}_v)+{\bf X}_c\boldsymbol\beta_c+{\bf D} \boldsymbol\alpha+\boldsymbol\varepsilon.   
\end{equation}
 
Given ${\bm \beta}_c$, and assuming that ${\boldsymbol\gamma}_v(\tau_t)$ has continuous derivatives of up to the second order, similar to the  LLDV discussed in the preceding section, we can derive the estimator   of ${\boldsymbol\gamma}_v(\tau_t)$ as follows: 
$$
\hat{\boldsymbol\gamma}_{v,{\bm \beta}_c}(\tau_0)={\bm{\Phi}}(\tau_0)({\bm Y}-{\bf X}_c\boldsymbol\beta_c), 
$$
where ${\bm{\Phi}}(\tau_0)=\left({\bf I}_{q+1}, {\bf 0 }_{(q+1) \times (q+1)}\right) \left({\bf M}^{\top}(\tau_0){\bf W}_h^{*}(\tau_0){\bf M}(\tau_0)\right)^{-1} {\bf M}^{\top}(\tau_0){\bf W}_h^{*}(\tau_0)$, and 
$$
{\bf M}(\tau_{0})=\left(\begin{array}{cc}
	{\bf z}_{v,11}^{\top} & \frac{\tau_{1}-\tau_{0}}{h} {\bf z}_{v,11}^{\top} \\
	\vdots & \vdots \\
	{\bf z}_{v,N1}^{\top} & \frac{\tau_{1}-\tau_{0}}{h} {\bf z}_{v,N1}^{\top} \\
	\vdots & \vdots \\
	{\bf z}_{v,NT}^{\top} & \frac{\tau_{T}-\tau_{0}}{h} {\bf z}_{v,NT}^{\top}\\
\end{array}\right).
$$
We refer to $\hat{\boldsymbol\gamma}_{v,{\bm \beta}_c}(\tau_0)$   as the PLLDV
estimator.  

Furthermore,  we define  
$$
{\bf S}(h)=\left(\begin{array}{cc}
	({\bf z}_{v,11}^\top, {\bf 0}_{1 \times (q+1)})[{\bf M}^{\top}(\tau_1) {\bf W}_h^{*}(\tau_1) {\bf M}(\tau_1)]^{-1} {\bf M}^{\top}(\tau_1) {\bf W}_h^{*}(\tau_1)\\
	\vdots \\
	({\bf z}_{v,N1}^\top, {\bf 0}_{1 \times (q+1)})[{\bf M}^{\top}(\tau_1) {\bf W}_h^{*}(\tau_1) {\bf M}(\tau_1)]^{-1} {\bf M}^{\top}(\tau_1) {\bf W}_h^{*}(\tau_1)\\
	\vdots \\
	({\bf z}_{v,NT}^\top, {\bf 0}_{1 \times (q+1)})[{\bf M}^{\top}(\tau_T) {\bf W}_h^{*}(\tau_T) {\bf M}(\tau_T)]^{-1} {\bf M}^{\top}(\tau_T) {\bf W}_h^{*}(\tau_T)\\
\end{array}\right).
$$
We can obtain the estimator of ${\bf B}({\bf Z}_v,{\bm \gamma}_v)$ as $\hat{\bf B}_{{\bm \beta}_c}({\bf Z}_v,{\bm \gamma}_v)={\bf S}(h)({\bm Y}-{\bf X}_c\boldsymbol\beta_c)$. 
By substituting $\hat{\bf B}_{{\bm \beta}_c}({\bf Z}_v,{\bm \gamma}_v)$ into Eq.(\ref{9}), we have 
\begin{equation} \label{10}
   \tilde {\bm Y}=\tilde {\bf X}_{c}{\bm \beta}_c+{\bf D}\bm\alpha+\boldsymbol\varepsilon,
\end{equation}
where $\tilde {\bm Y}=({\bf I}_{NT}-{\bf S}(h)){\bm Y}$, $\tilde {\bf X}_c=({\bf I}_{NT}-{\bf S}(h)){\bf X}_c$. The profile least square estimator of $\bm \alpha$  is given by
\begin{equation} \label{11}
  \hat{\bm \alpha}_{\beta_c}=({\bf D}^\top{\bf D})^{-1}{\bf D}^\top(\tilde {\bm Y}-\tilde {\bf X}_c \bm\beta_c).
\end{equation}
Denote ${\bf P}_{\bf D}={\bf D}({\bf D}^\top {\bf D})^{-1} {\bf D}^\top$. Replacing ${\bm \alpha}$ in   Eq.(\ref{10}) by $\hat{\bm \alpha}_{{\bm \beta}_c}$,  we can rewrite Eq.(\ref{10})  as
\begin{equation}  \label{12}
    \bar {\bm Y}=\bar {\bf X}_{c}{\bm \beta}_c+\boldsymbol\varepsilon,
\end{equation}
where $\bar {\bm Y}=({\bf I}_{NT}-{\bf P}_{\bf D})\tilde {\bf Y}$ and $\bar {\bf X}_c=({\bf I}_{NT}-{\bf P}_{\bf D})\tilde {\bf X}_c$.  Consequently, the least squares estimator of   ${\bm \beta}_c$   is obtained as 
$$ 
\hat{\bm \beta}_c=(\bar{\bf X}_c^\top \bar{\bf X}_c)^{-1}\bar{\bf X}_c^\top \bar {\bm Y}.
$$

By substituting $\hat{\bm \beta}_c$ into $\hat{\bm \alpha}_{\beta_c}$, $\hat{\boldsymbol\gamma}_{{\bm \beta}_c}(\tau_0)$ and $\hat{\bf B}_{{\bm \beta}_c}({\bf Z}_v,{\bm \gamma}_v)$, the final estimators of $\bm\alpha$, ${\boldsymbol\gamma}(\tau_0)$ and ${\bf B}({\bf Z}_v,{\bm \gamma}_v)$ are respectively given as
$$
\hat{\bm \alpha}=({\bf D}^\top{\bf D})^{-1}{\bf D}^\top(\tilde {\bm Y}-\tilde {\bf X}_c \hat{\bm\beta}_c),
$$
$$
\hat{\boldsymbol\gamma}_v(\tau_0)={\bm{\Phi}}(\tau_0)({\bm Y}-{\bf X}_c \hat{ \boldsymbol\beta}_c), 
$$
and 
$$
\hat{\bf B}({\bf Z}_v,{\bm \gamma}_v)={\bf S}(h)({\bm Y}-{\bf X}_c \hat{ \boldsymbol\beta}_c).
$$
The residual vector and the residual sum of squares are provided respectively as follows:
$$
\hat{\boldsymbol\varepsilon}_{\bf PL}(h)=({\bf I}_{NT}-{\bf L}(h))({\bf I}_{NT}-{\bf P}_{\bf D})({\bf I}_{NT}-{\bf S}(h)){\bm Y}, 
$$
and 
\begin{equation} \label{13}
    { RSS}_{\bf PL}(h)=\hat{\boldsymbol\varepsilon}_{\bf PL}^\top(h)\hat{\boldsymbol\varepsilon}_{\bf PL}(h)=\bar{\bm Y}^\top({\bf I}_{NT}-{\bf L}(h))^\top({\bf I}_{NT}-{\bf L}(h))\bar{\bm Y}, 
\end{equation}
where ${\bf L}(h)=\bar {\bf X}_c(\bar {\bf X}_c^\top\bar {\bf X}_c)^{-1}\bar {\bf X}_c^\top$.

\section{Asymptotic properties}
  Throughout this analysis and in subsequent discussions, we define  $\mu_i=\int u^i K(u) du$ and $\nu_i=\int u^i K^2(u) du$ for $i=0,1,2$. 
To derive the  asymptotic normality properties of the 2SLS-PLLDV estimators, $\hat{\bm\beta}_c$ and $\hat{\bm\gamma}_v(\cdot)$, we outline the subsequent regularity assumptions for   Model (\ref{2}).

\noindent\textbf{Assumption 1} ~ The kernel function $K(\cdot)$ is a bounded symmetric and Lipschitz continuous with a compact support $[-1, 1]$. Typically, $K(\cdot)$ satisfies that
$$
\int K(u) du = 1, \quad \int u K(u) du = 0, \quad \int u^2 K(u) du < \infty.
$$ 

\noindent\textbf{Assumption 2} ~The covariates $\{{\bf x}_{c,it}, {\bf x}_{v,it}\}$ satisfy certain regularity conditions: (i) $\{({\bf X}_{c,i}, {\bf X}_{v,i}), $  $ i \geq 1\}$ is a sequence of independent and identically distributed (i.i.d.) variables, where ${\bf X}_{c,i}= ({\bf x}_{c,it}, t \geq 1)$ and ${\bf X}_{v,i} = ({\bf x}_{v,it}, t \geq 1)$. Furthermore, for every $i \geq 1$, $\{({\bf x}_{c,it}, {\bf x}_{v,it}), t \geq 1\}$ is stationary and $\alpha$-mixing with mixing coefficient $\alpha_k$ satisfying $\alpha_k = O(k^{-\tau_0})$, where $\tau_0 > \frac{\delta + 2}{\delta}$ for some $\delta > 0$ involved in (ii) below. (ii) $E({\bf x}_{c,it})={\bf 0}_{p-q}$, $E({\bf x}_{v,it})={\bf 0}_q$, and there exist positive definite matrices $\bm\Sigma_{{\bf x}_c} = E({\bf x}_{c,it}{\bf x}_{c,it}^\top)$, $\bm\Sigma_{{\bf z}_v} = E({\bf z}_{v,it}{\bf z}_{v,it}^\top)$, $\bm\Sigma_{{\bf z}_v{\bf x}_c} = E({\bf z}_{v,it}{\bf x}_{c,it}^\top)$, and $\bm\Sigma=\bm\Sigma_{{\bf x}_c}-\bm\Sigma_{{\bf z}_v {\bf x}_c}^\top \bm\Sigma_{{\bf z}_v}^{-1}\bm\Sigma_{{\bf z}_v {\bf x}_c}$. (iii) $E\|{\bf x}_{c,it}\|^{2(2+\delta)} < \infty$, $E\|{\bf z}_{v,it}\|^{2(2+\delta)} < \infty$ for some $\delta > 0$, where $\|\cdot\|$ is the $L_2$-distance. 

\noindent\textbf{Assumption 3}~ The error process $\{\varepsilon_{it}\}$ is independent of both ${\bf x}_{c,it}$ and $ {\bf x}_{v,it}$ with $E(\varepsilon_{it}) = 0$, $E(\varepsilon_{it}^2) = \sigma_{\bm \varepsilon}^2 < \infty$, $E(|\varepsilon_{it}|^{2+\delta}) < \infty$ for some $\delta>0$. For every $i$, $\bm \varepsilon_{i}$s are stationary and $\alpha$-mixing sequence,  as do the pairs  $({\bf X}_{c,i}, {\bf X}_{v,i})$, where $\bm \varepsilon_{i}=(\varepsilon_{it},t \geq 1)$. Furthermore, $c_{\bm \varepsilon}(t) = E(\varepsilon_{is} \varepsilon_{i,s+t})>0$, and ${\bf C}_{{\bf z}_v}(t) = E({\bf z}_{v,is} {\bf z}_{v,i,s+t}^\top)$, ${\bf C}_{{\bf x}_c}(t) = E({\bf x}_{c,is} {\bf x}_{c,i,s+t}^\top)$, ${\bf C}_{{\bf z}_v{\bf x}_c}(t) = E({\bf z}_{v,is}  {\bf x}_{c,i,s+t}^\top)$ are positive definite matrices. 

\noindent\textbf{Assumption 4}~ ${\bm \gamma}_v(\cdot)$ has continuous derivatives up to the second order.

\noindent\textbf{Assumption 5} ~ As $N,T \to \infty$ and $h \to 0$, the bandwidth $h$ satisfies that $Th \to \infty$. 

\noindent\textbf{Assumption 6} ~ The spatial adjacency matrix ${\bf W}_N$ is nonstochastic and uniformly bounded (UB) in both row and column sums in absolute value.
\\ \\
\textbf{Theorem 3.1.}  Under Assumptions 1-6, as $N,T \to \infty$, we have 
$$
\sqrt{NT}(\hat{\bm{\beta}}_c - \bm{\beta}_c) \overset{d}{\rightarrow} N({\bf 0}_{p-q}, \bm\Sigma^{-1}\bm\Omega_{\bm\varepsilon}\bm\Sigma^{-1}),
$$
where $\bm\Sigma=\bm\Sigma_{{\bf x}_c}-\bm\Sigma_{{\bf z}_v {\bf x}_c}^\top \bm\Sigma_{{\bf z}_v}^{-1}\bm\Sigma_{{\bf z}_v {\bf x}_c}$, 
$\bm\Omega_{\bm\varepsilon}=\sum\limits_{t=-\infty}^\infty c_{\bm\varepsilon}(t)\left({\bf C}_{{\bf x}_c}(t)-{\bf C}_{{{\bf x}_c}{{\bf z}_v}}(t)\bm\Sigma_{{\bf z}_v}^{-1} \bm\Sigma_{{\bf z}_v {\bf x}_c}\right.\\
\left.-\bm\Sigma^\top_{{\bf z}_v {\bf x}_c}\bm\Sigma_{{\bf z}_v}^{-1}{\bf C}_{{{\bf z}_v}{{\bf x}_c}}(t)+\bm\Sigma^\top_{{\bf z}_v {\bf x}_c}\bm\Sigma_{{\bf z}_v}^{-1}{\bf C}_{{\bf z}_v}(t)\bm\Sigma_{{\bf z}_v}^{-1} \bm\Sigma_{{\bf z}_v {\bf x}_c}\right)$.  
 \\ 
 \\
\textbf{Theorem 3.2.} Under Assumptions 1-6, for given $0 < \tau_0 < 1$, as $N,T \to \infty$, we have
$$
\sqrt{NTh}[\hat{\bm\gamma}_v(\tau_0)-{\bm\gamma}_v(\tau_0)- b(\tau_0)h^2 + o_p(h^2)] \xrightarrow{d} N(\mathbf{0}_{q+1}, \bm\Sigma_{{\bf z}_v}^{-1} \bm\Omega \bm\Sigma_{{\bf z}_v}^{-1}),
$$
  where $b(\tau_0) = \frac{1}{2} \mu_2 {\bm\gamma}_v''(\tau_0)$, ${\bf \Omega}=\nu_0 \sum\limits_{t=-\infty}^{\infty}{\bf C}_{{\bf z}_v}(t)c_{\bm \varepsilon}(t)$.

The proofs of Theorems 3.1 and 3.2 are 
given in Appendix.

\section{Test for time-varying regression coefficients}
In this section, the test statistics  are constructed based on the residual sums of squares in order to address the testing problems (\ref{3}).

\subsection{Construction of the test statistic}
If the alternative hypothesis $H_{1}$  holds true, Model (1) becomes a special case of Model (2) when  $q = p$. Consequently, we can derive   the residual vector and the residual sum of squares of Model (1), respectively,  given by
\begin{equation} \label{14}
    \hat{\boldsymbol\varepsilon}_{\bf TV}(h)=({\bf I}_{NT}-{\bf P}_{\bf D})({\bf I}_{NT}-{\bf S}(h)){\bm Y},
\end{equation}
and 
\begin{equation} \label{15}
   { RSS}_{\bf TV}(h)=\hat{\boldsymbol\varepsilon}_{\bf TV}^\top(h)\hat{\boldsymbol\varepsilon}_{\bf TV}(h)=\bar{\bm Y}^\top\bar{\bm Y}.
\end{equation}
Utilizing ${ RSS}_{\bf TV}(h)$ in Eq.(\ref{15}) and ${ RSS}_{\bf PL}(h)$ in Eq.(\ref{13}),  we construct the test statistic similar to that of Mei and Chen (2022). The test statistic is defined as
\begin{equation} \label{16}
    W=\frac{NT}{2}\frac{{ RSS}_{\bf PL}(h)-{RSS}_{\bf TV}(h)}{{ RSS}_{\bf TV}(h)}.
\end{equation}
Clearly, a large value of the test statistic $W$ indicates that the null hypothesis should be rejected. Let $w$  denote the observed value of $W$, and the p-value is computed as
\begin{equation} \label{17}
   {\rm p}=P_{H_0}(W \geq w),
\end{equation}
where $P_{H_0}(\cdot)$  represents the probability computed under the null hypothesis $H_0$. For a given significance level $\alpha$, we reject the null hypothesis $H_0$ if ${\rm p}<\alpha$, otherwise not reject $H_0$.

\subsection{Bootstrap method based on residual to calculate p-value}
\par
 To calculate the  p-value of the test statistic, a commonly adopted  method is to establish the asymptotic null distribution. However, deriving this distribution is challenging due to the presence of spatial lag terms. In the field of non-parametric and semi-parametric regression, researchers have utilized the bootstrap method based on residuals to compute   p-values. For instance, Tian LL {\it et al.} (2024), Li and Mei (2013, 2016) and Wei {\it et al.} (2017) have proposed residual-based bootstrap procedures to approximate the   p-value for spatial autoregressive models. Similarly, in  Model (2), we also employ the residual-based bootstrap method to compute the \( p \)-value of our test statistic. 
 
 The steps involved are outlined below:
\\
{\bf Step 1.} Based on the data set $\{y_{it}, {\bf x}_{it}\},$  $i=1,2,\cdots,N$; $t=1,2,\cdots,T$, compute the observed value $w$ of statistics $W$ by Eq.(\ref{16}). Then under alternative hypotheses $H_1$, yield the estimators of the residuals $\hat{\bm \varepsilon}_{\bf TV}(h)=(\hat\varepsilon_{11},\cdots,\hat\varepsilon_{N1},\cdots,\hat\varepsilon_{NT})^\top$ by Eq. (\ref{14}) and centralize it to obtain
$\hat{\boldsymbol\varepsilon}_c=(\hat\varepsilon_{11}-\bar{\hat{\varepsilon}},\hat\varepsilon_{12}-
\bar{\hat{\varepsilon}},\cdots,\hat\varepsilon_{NT}-\bar{\hat{\varepsilon}})^\top$, where
$\bar{\hat{\varepsilon}}=\frac{1}{NT}\sum\limits_{i=1}^N\sum\limits_{T=1}^T\hat\varepsilon_{it}$.
\\
{\bf Step 2.} Draw a bootstrap residuals
$\hat{\boldsymbol\varepsilon}^*=(\hat\varepsilon_{11}^*,\hat\varepsilon_{12}^*,\cdots,\hat\varepsilon_{NT}^*)$
 from the empirical distribution
function of $\hat{\boldsymbol\varepsilon}_c$.
\\
{\bf Step 3.} Generate ${\bm Y}^*=({\bf I}_{NT}-\hat{\bm\rho}_{NT} {\bf W})^{-1}(\hat{\bf B}({\bf X}_v,{\bm \beta}_v)+{\bf X}_c\hat{\bm \beta}_c +\hat{\boldsymbol\varepsilon}^*)$,
and calculate the bootstrap test statistics $W^*$ by Eq.(\ref{16}), we get
\begin{equation}  \label{18}
    W^*=\frac{NT}{2}\frac{{ RSS}_{\bf PL}^*(h)-{RSS}_{{\bf TV}}^*(h)}{{ RSS}_{{\bf TV}}^*(h)},
\end{equation}
where ${RSS}_{{\bf TV}}^*(h)$ is the estimators by re-fitting the model (1) based on the bootstrap sample $\left\{{\bm Y}^*,{\bf X}\right\}$, and ${ RSS}_{\bf PL}^*(h)$ is the estimators by re-fitting the model (4) based on the bootstrap sample $\left\{{\bm Y}^*,{\bf X}_c,{\bf X}_v\right\}$;
\\
{\bf Step 4.} Repeat the step 2 and step 3  a total of $k$ times and yield 
bootstrap values of the test statistics $W$, denoted as  $W_{1}^*, W_{2}^*,\cdots, W_{k}^*$. Subsequently,  the  p-value of the test is derived by
$$
\hat {\rm p}=\frac{\sharp \{W_{j}^*: W_{j}^* \geq w \}}{k},~~j=1,2,\cdots,k,
$$
where $w$ is the observed value of the test statistic $W$ obtained in step 1 and $\sharp A$ denotes the number of the elements in set $A$.

\section{Simulation studies}
In this section, we  assess the performance of the proposed 2SLS-PLLDV estimation and testing framework in detecting linear relationships between explanatory variables and the response variable, particularly in the context of time-varying regression coefficients. This evaluation is conducted through a series of simulation studies.

We define a grid consisting $m$ rows and $m$ columns, resulting in $N=m^2$ location observations  distributed across   the $m\times m$ regions. In the simulations, the Gaussian kernel function $K(\frac{z}{h})=\frac{1}{\sqrt{2\pi}}\exp(-\frac{z^2}{2h^2})$ with the ROT (the Rule of Thumb) bandwidth $h=s_{\tau}(NT)^{-1/5}$ is used to generate the kernel weights, where $s_{\tau}$ is the sample standard deviation of a vector $\tau=(\tau_1,\tau_2,\cdots,\tau_T)^\top$. 

\subsection{Simulation of estimation method}
\par
We generate the experimental data from the following model:
\begin{equation} \label{19}
    \begin{array}{lll}
	y_{i t}=&\rho(\tau_t)\sum\limits_{j \neq i} w_{i j} y_{j t}+\beta_1(\tau_{t})x_{it1}+\beta_2(\tau_{t})x_{it2}+\beta_3x_{it3}\\
	&+\beta_4x_{it4}+\alpha_{i}+\varepsilon_{i t},i=1,2, \cdots, N; t=1,2, \cdots,T.
\end{array}  
\end{equation}
where the $x_{itj},j=1,2,3,4$ and the fixed effects $\alpha_i$ are  drawn from independent standard normal distributions and uniform distributions. Specifically, $x_{it1}=1$, $x_{it2}\sim N(0, 1)$, $x_{it3}\sim N(0, 1)$, $x_{it4}\sim N(0, 1)$ and $\alpha_i \sim U(0,1)$. The relationship coefficients are given by 
$$
\beta_1(\tau_{t})=4\tau_t, ~~\beta_2(\tau_{t})=(\tau_t+1)^2,~~
    \beta_3=-5, ~~
    \beta_4=5. 
$$

To  show  the performance of various forms of the spatial adjacency matrix, spatial lag coefficient, and error term distributions, we consider two types of row-standardized spatial adjacency matrices: the Rook and Queen adjacency matrices. Additionally, we consider the following two forms of the spatial lag coefficient and three error distribution types, all with a mean of  $\mu=0$ and a variance $\sigma^2=1$:
$$
\rho_1(\tau_{t})=-0.6\sin^2(2\pi \tau_t),\quad 
\rho_2(\tau_{t})=0.6\sin^2(2\pi \tau_t),
$$

$$
\left\{
\begin{array}{l}   \varepsilon_{it} \sim N(0,1)\\
\varepsilon_{it} \sim
U(-\sqrt{3},\sqrt{3}), \\
\varepsilon_{it} \sim
\frac{1}{2}\chi^2(2)-1\\
\end{array}
\right.
$$
where the range of $\rho_1(\tau_{t})$ and $\rho_2(\tau_{t})$ are $[-0.6,0]$ and $[0,0.6]$, respectively. We use the average bias, standard deviation (SD) as the measure for estimation accuracy of constant coefficients $\beta_3$ and $\beta_4$, and the average mean squared errors (AMSE) of the estimates of the time-varying coefficients to assess their performance. The AMSE of $\hat\rho(\cdot)$ and $\hat{\bm\beta}(\cdot)$ are defined as 
$$
AMSE(\hat{\rho}(\cdot)) = \frac{1}{T} \sum_{t=1}^{T} \left\{ \frac{1}{n_{sim}} \sum_{i=1}^{n_{sim}} \left[ \hat{\rho}^{(i)} (\tau_t)-\rho(\tau_t)\right]^2 \right\},
$$
$$
AMSE(\hat{\beta}_j(\cdot)) = \frac{1}{T} \sum_{t=1}^{T} \left\{ \frac{1}{n_{sim}} \sum_{i=1}^{n_{sim}} \left[ \hat{\beta}_j^{(i)} (\tau_t) - \beta_j (\tau_t)\right]^2 \right\}, \quad j = 1, 2,
$$
where $\hat\rho^{(i)}(\cdot)$ and $\hat{\bm\beta}^{(i)}(\cdot)$ denote the estimates of $\hat\rho(\cdot)$ and ${\bm\beta}(\cdot)$ in the $i$-th replication, with $n_{sim}$ denoting the number of replications, chosen as   $n_{sim}=500$. Additionally, four distinct sample sizes   $N  T$, specifically   $N=10^2$ and $12^2$, with $T=5$ and $10$. The simulation results   for the  AMSE  of the time-varying coefficients $\rho(\cdot)$ and $\bm\beta(\cdot)$ are summarized in Table 1.  Furthermore, Table 2 contains the average bias, standard deviation (SD) of the constant coefficients $\beta_3$ and $\beta_4$. Finally, we present the average estimates curves of $\rho(\tau_t)$, $\beta_1(\tau_t)$ and $\beta_2(\tau_t)$ under the setting of $\rho_1(\tau_t)$, Rook adjacency matrix and $\varepsilon_{it} \sim N(0,1)$ for $N=144$ and $T=10$ in Figure 1, alongside their corresponding true curves for comparison.  

\begin{table}
	\centering
	\scriptsize
	\vskip -0.3cm
	\caption{AMSE for estimators of $\rho(\tau_t)$, $\beta_1(\tau_t)$ and $\beta_2(\tau_t)$ under Model \eqref{19}}
	\label{table1}
    \renewcommand{\arraystretch}{0.5}  
	\resizebox{\textwidth}{!}{  
		\begin{tabular}{ccccccc}
			\toprule
			$\rho(\tau_t)$ & Adjacency matrix& Distribution of $\varepsilon_{i t}$ & $(T,N)$ & $\hat{\rho}(\tau_t)$ & $\hat{\beta}_1(\tau_t) $ & $\hat{\beta}_2( \tau_t$) \\
			\midrule
		\multirow{24}{*}{$\rho_1(\tau_t)$} & \multirow{12}{*}{Rook} & \multirow{4}{*}{$N(0,1)$} & (5,100) &   0.0117 & 0.0319 & 0.0192 \\
		& & & (5,144) & 0.0113 & 0.0315 & 0.0156 \\
		& & & (10,100) & 0.0062 & 0.0118 & 0.0095  \\
		& & & (10,144) & 0.0052 & 0.0087 & 0.0072 \\
		& & \multirow{4}{*}{$U(\sqrt{3},\sqrt{3})$} & (5,100) &0.0118 & 0.0314 & 0.0198 \\
		& & & (5,144) & 0.0114 & 0.0335 & 0.0144 \\
		& & & (10,100) & 0.0061 & 0.0109 & 0.0099\\
		& & & (10,144) & 0.0052 & 0.0083 & 0.0070\\
		& & \multirow{4}{*}{$\frac{1}{2} \chi^{2}(2)-1$} & (5,100) & 0.0118 &0.0328  &0.0197  \\
		& & & (5,144) & 0.0111 &0.0328 & 0.0143 \\
		& & & (10,100) & 0.0062 & 0.0114 & 0.0094 \\
		& & & (10,144) &0.0052 &0.0083 & 0.0069 \\
		& \multirow{12}{*}{Queen} & \multirow{4}{*}{$N(0,1)$} & (5,100) &0.0206  & 0.0672 &  0.0148\\
		& & & (5,144) & 0.0199 &0.0669  &0.0117 \\
		& & & (10,100) &0.0092& 0.0113 & 0.0083 \\
		& & & (10,144) & 0.0080&0.0088  &0.0061  \\
		& & \multirow{4}{*}{$U(\sqrt{3},\sqrt{3})$} & (5,100) & 0.0209 & 0.0652 & 0.0156  \\
		& & & (5,144) & 0.0202 & 0.0689 &0.0107 \\
		& & & (10,100) & 0.0093 &0.0112 & 0.0088\\
		& & & (10,144) & 0.0080 & 0.0087 & 0.0060 \\
		& & \multirow{4}{*}{$\frac{1}{2} \chi^{2}(2)-1$} & (5,100) &0.0203 & 0.0645& 0.0155 \\
		& & & (5,144) & 0.0197 & 0.0679 & 0.0103 \\
		& & & (10,100) & 0.0092 & 0.0110 &0.0081 \\
		& & & (10,144) & 0.0079 & 0.0087 & 0.0059 \\
   \multirow{24}{*}{$\rho_2(\tau_t)$} & \multirow{12}{*}{Rook} & \multirow{4}{*}{$N(0,1)$} & (5,100) & 0.0059 &0.0710  &0.0179 \\
		& & & (5,144) & 0.0052 & 0.0637 & 0.0145 \\
		& & & (10,100) & 0.0061 & 0.0234 &0.0097 \\
		& & & (10,144) & 0.0051 & 0.0202 & 0.0074 \\
		& & \multirow{4}{*}{$U(\sqrt{3},\sqrt{3})$} & (5,100) & 0.0060 &0.0714 &0.0194  \\
		& & & (5,144) & 0.0053 & 0.0645 &0.0133\\
		& & & (10,100) & 0.0062 & 0.0222 & 0.0100 \\
		& & & (10,144) & 0.0051 & 0.0195& 0.0070 \\
		& & \multirow{4}{*}{$\frac{1}{2} \chi^{2}(2)-1$} & (5,100) & 0.0059&0.0727 &0.0192  \\
		& & & (5,144) & 0.0051 & 0.0633 & 0.0131 \\
		& & & (10,100) & 0.0061 & 0.0227 & 0.0095 \\
		& & & (10,144) & 0.0050 & 0.0193 & 0.0069 \\
		& \multirow{12}{*}{Queen} & \multirow{4}{*}{$N(0,1)$} & (5,100) & 0.0070&0.1300 &0.0161  \\
		& & & (5,144) & 0.0063 & 0.1135 &0.0124 \\
		& & & (10,100) & 0.0069 & 0.0348& 0.0089 \\
		& & & (10,144) & 0.0056 & 0.0284 & 0.0066 \\
		& & \multirow{4}{*}{$U(\sqrt{3},\sqrt{3})$} & (5,100) & 0.0070 & 0.1283 & 0.0173 \\
		& & & (5,144) & 0.0064 & 0.1161& 0.0120 \\
		& & & (10,100) & 0.0069 & 0.0340 & 0.0092 \\
		& & & (10,144) & 0.0056 & 0.0282& 0.0064 \\
		& & \multirow{4}{*}{$\frac{1}{2} \chi^{2}(2)-1$} & (5,100) & 0.0071 & 0.1307 & 0.0171\\
		& & & (5,144) & 0.0063 & 0.1164 &0.0113 \\
		& & & (10,100) & 0.0069 & 0.0362 & 0.0087 \\
		& & & (10,144) & 0.0054 & 0.0278 & 0.0062  \\
			\bottomrule
		\end{tabular}
}
\end{table}

\begin{table}
	\centering
	\scriptsize
	\vskip -0.3cm
	\caption{Bias and SD for estimators of $\beta_3$ and $\beta_4$ under Model \eqref{19}} 
	\label{table2}
    \renewcommand{\arraystretch}{0.5}  
	\resizebox{\textwidth}{!}{  
		\begin{tabular}{ccccccccc}
			\toprule
			& & & & \multicolumn{2}{c}{$\hat{\beta}_3$} & \multicolumn{2}{c}{$\hat{\beta}_4$} \\
  \cmidrule{5-6}  \cmidrule{7-8}
			$\rho(\tau_t)$  &  Adjacency matrix & Distribution of $\varepsilon_{i t}$ & $(T,N)$ & Bias & SD & Bias & SD \\
			\midrule
  \multirow{24}{*}{$\rho_1(\tau_t)$} & \multirow{12}{*}{Rook} & \multirow{4}{*}{$N(0,1)$} & (5,100) & -0.0328& 0.0652 & 0.0322 & 0.0615\\
		& & & (5,144) & -0.0317 &0.0516  & 0.0310 & 0.0526   \\
		& & & (10,100) &-0.0275 &0.0372 &0.0255 & 0.0363 \\
		& & & (10,144) & -0.0233 &0.0303 &  0.0209 & 0.0315\\
		& & \multirow{4}{*}{$U(\sqrt{3},\sqrt{3})$} & (5,100) & -0.0320 & 0.0598 &  0.0365 & 0.0624    \\
		& & & (5,144) & -0.0297 & 0.0531 & 0.0297 & 0.0518 \\
		& & & (10,100) & -0.0265 &0.0371 & 0.0265 & 0.0376  \\
		& & & (10,144) & -0.0277 & 0.0321& 0.0266 & 0.0306 \\
		& & \multirow{4}{*}{$\frac{1}{2} \chi^{2}(2)-1$} & (5,100) & -0.0317 & 0.0644 & 0.0331 &  0.0652 \\
		& & & (5,144) &-0.0345 & 0.0527 &  0.0332 & 0.0529   \\
		& & & (10,100) &-0.0214 &0.0372 &  0.0258 &  0.0365   \\
		& & & (10,144) & -0.0222 & 0.0328 & 0.0214 & 0.0324 \\
          
   & \multirow{12}{*}{Queen} & \multirow{4}{*}{$N(0,1)$} & 
   (5,100) & -0.0187 &  0.0579   & 0.0182 & 0.0559 \\
		& & & (5,144) & -0.0169 & 0.0469 &0.0148 &0.0503 \\
		& & & (10,100) &-0.0173 &  0.0353 &0.0141 & 0.0346 \\
		& & & (10,144) & -0.0138  & 0.0278 &0.0116 & 0.0285 \\
        
		& & \multirow{4}{*}{$U(\sqrt{3},\sqrt{3})$} 
        & (5,100) & -0.0167  & 0.0607 & 0.0214 & 0.0556 \\
		& & & (5,144) & -0.0191  & 0.0481  & 0.0133 & 0.0483 \\
		& & & (10,100) &  -0.0110 &  0.0352  & 0.0164 & 0.0359\\
		& & & (10,144) & -0.0138  & 0.0299  & 0.0144 & 0.0287 \\
       
		& & \multirow{4}{*}{$\frac{1}{2} \chi^{2}(2)-1$} & (5,100) & 0.0331 & 0.0652 & 0.0197 &0.0580 \\
		& & & (5,144) & -0.0167  & 0.0607  & 0.0203 &0.0465\\
		& & & (10,100) & 0.0258 & 0.0365 &0.0158&0.0346  \\
		& & & (10,144) & 0.0214&0.0324 & 0.0124 & 0.0305 \\
    
   \multirow{24}{*}{$\rho_2(\tau_t)$}   & \multirow{12}{*}{Rook} 
  & \multirow{4}{*}{$N(0,1)$} & (5,100) & -0.0308 & 0.0655 &0.0320   &  0.0583  \\
		& & & (5,144) & -0.0282 & 0.0500& 0.0257   & 0.0527 \\
		& & & (10,100) & -0.0269 & 0.0373& 0.0211   & 0.0356 \\
		& & & (10,144) & -0.0198 & 0.0305& 0.0186   & 0.0302 \\
		& & \multirow{4}{*}{$U(\sqrt{3},\sqrt{3})$} & (5,100) & -0.0284 & 0.0638 & 0.0366   & 0.0610 \\
		& & & (5,144)  & -0.0275  & 0.0511 & 0.0253   & 0.0496 \\
		& & & (10,100) & -0.0250 & 0.0366  & 0.0245  &  0.0372  \\
		& & & (10,144) & -0.0200 & 0.0313 & 0.0210  &  0.0305 \\
		& & \multirow{4}{*}{$\frac{1}{2} \chi^{2}(2)-1$} & (5,100) & -0.0272 & 0.0618 & 0.0269  &  0.0625\\
		& & & (5,144)& -0.0289 & 0.0528 &  0.0299   & 0.0501\\
		& & & (10,100) & -0.0186 & 0.0368&0.0237  &  0.0360  \\
		& & & (10,144) & -0.0207 & 0.0322 & 0.0191  &  0.0310 \\

  & \multirow{12}{*}{Queen} & \multirow{4}{*}{$N(0,1)$} &(5,100)&  -0.0172   & 0.0616 &0.0201&0.0547\\
		& & & (5,144) & -0.0147  &  0.0483  & 0.0123& 0.0487\\
		& & & (10,100) & -0.0166   & 0.0358 & 0.0114& 0.0346 \\
		& & & (10,144) &  -0.0112   & 0.0293  & 0.0101 & 0.0293\\
		& & \multirow{4}{*}{$U(\sqrt{3},\sqrt{3})$} & (5,100) & -0.0152   & 0.0592 &0.0222 &0.0568 \\
		& & & (5,144) &-0.0138   & 0.0488 & 0.0125&0.0465\\
		& & & (10,100) &  -0.0139  &  0.0356  &0.0135 & 0.0364 \\
		& & & (10,144) &-0.0112  &  0.0299 & 0.0124 & 0.0294 \\
		& & \multirow{4}{*}{$\frac{1}{2} \chi^{2}(2)-1$} & (5,100) & -0.0151    &0.0584 & 0.0156 & 0.0598\\
		& & & (5,144) & -0.0144   & 0.0498 & 0.0160 & 0.0471 \\
		& & & (10,100) &-0.0078   & 0.0358 & 0.0129&0.0350 \\
		& & & (10,144) & -0.0114   & 0.0303 & 0.0105 & 0.0298 \\
		\bottomrule
		\end{tabular}
	}
\end{table}

\begin{figure}[htbp]
    \centering
    \begin{subfigure}{0.5\textwidth}
        \centering
        \includegraphics[width=\textwidth]{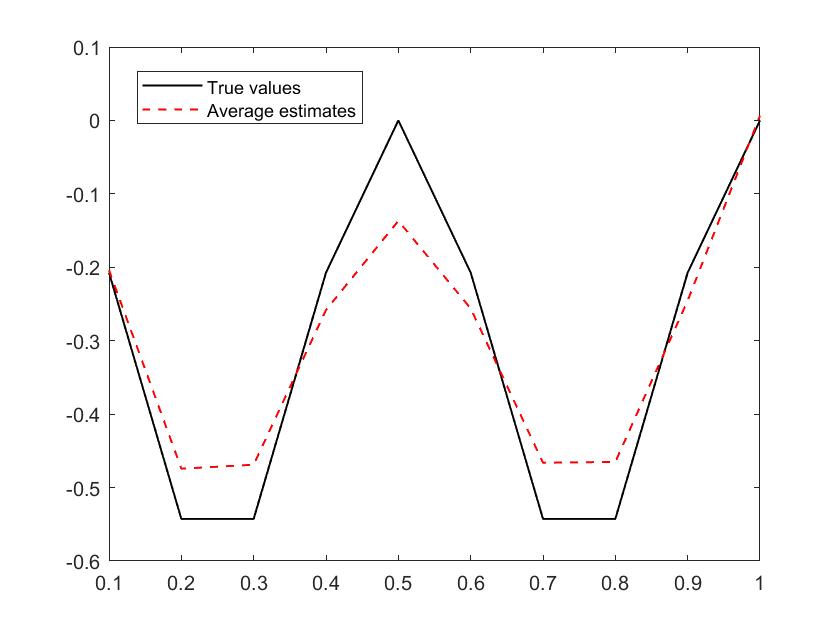}
        \caption{$\rho_1(\tau_t)$}
        \label{fig:subfig1}
    \end{subfigure}
    \begin{subfigure}{0.5\textwidth}
        \centering
        \includegraphics[width=\textwidth]{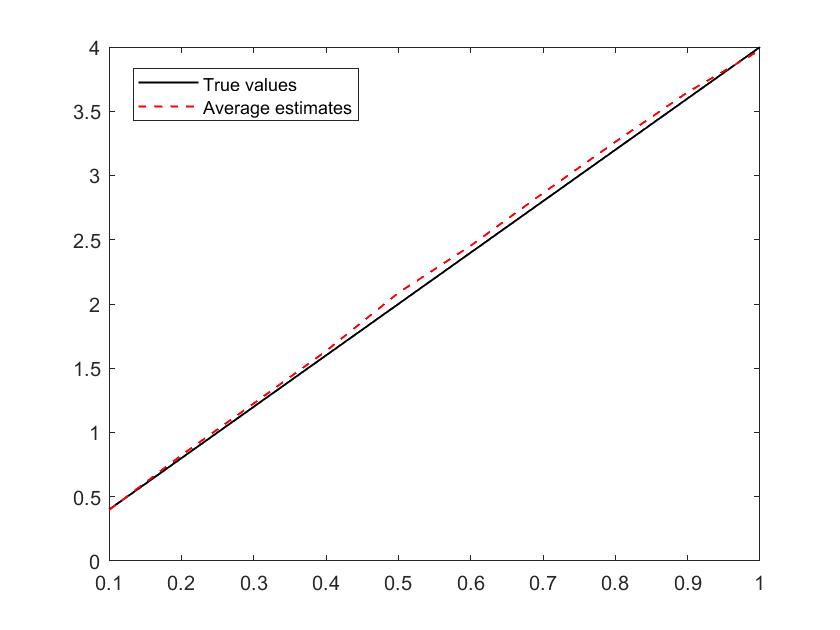}
        \caption{$\beta_1(\tau_t)$}
        \label{fig:subfig2}
    \end{subfigure}
    \begin{subfigure}{0.5\textwidth}
        \centering
        \includegraphics[width=\textwidth]{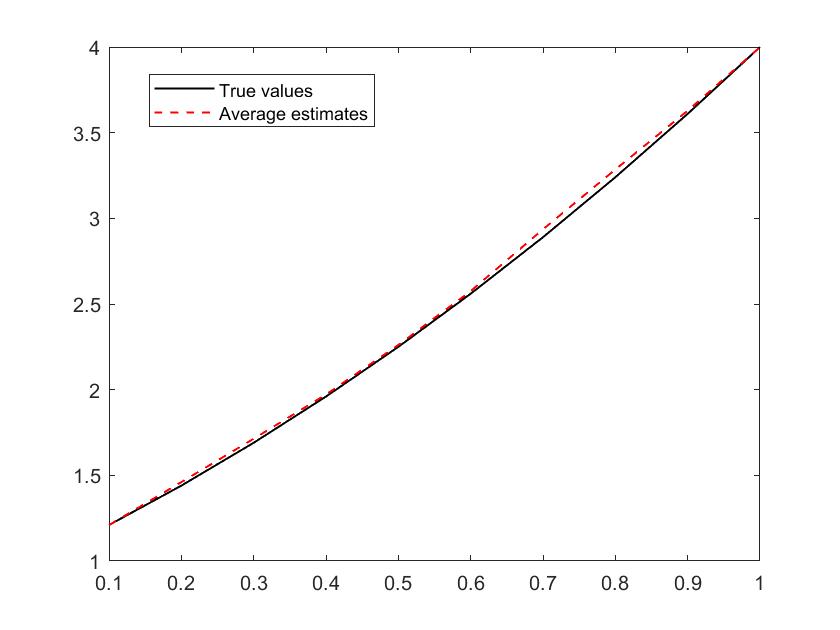}
        \caption{$\beta_2(\tau_t)$}
        \label{fig:subfig3}
    \end{subfigure}
    \caption{The average estimates curves and true  curves of $\rho(\tau_t)$, $\beta_1(\tau_t)$ and $\beta_2(\tau_t)$}
    \label{fig:threesubfigures}
\end{figure}

From Table 1 and Table 2, we can observe that the AMSE for estimators of the time-varying coefficients $\rho(\tau_t)$ and $\bm\beta(\tau_t)$ is low, and similarly,  the Bias and SD for estimators of the constant coefficients $\beta_3$ and $\beta_4$ are small. Notably, as the sample size dimensions either $N$ or $T$ increases, the values of AMSE and SD decrease.  Across  the three distributions of the error term,  no significant differences are apparent in the measures' values. Furthermore, the results remain consistent regardless of the forms of the spatial lag coefficients and spatial adjacency matrices, highlighting the effectiveness and robustness of our proposed estimation method.

\subsection{Simulation on  the proposed test}
In this simulation, the data generation model retains the similar structure as that presented in Eq. (\ref{19}), where the time-varying coefficients $\beta_1(\tau_t)$ and $\beta_2(\tau_t)$  are maintained, while the constant coefficients $\beta_3$ and $\beta_4$  are modified as 
$$
\left\{
\begin{array}{l}
\beta_3(\tau_{t})=-5+ce^{\tau_t}, \\
\beta_4(\tau_{t})=5+ c\sin(\pi\tau_t), \\
\end{array}
\right.
$$
where $c$  takes values within the interval [0,1],  characterizing the degree  of deviation from the null hypothesis, with larger values of $c$ indicating greater deviation.  The remaining experimental settings, including the distributions of the independent variables, fixed effects and error term, the functional forms of the time-varying coefficients associated with  the spatial lag term and the types of the spatial adjacency matrices, are  the same as  those used in the simulation outlined in Section 5.1. 

For each experimental setting,  we conduct  $n_{sim}=500$ simulation runs with varying sample sizes $NT$ (where $N=8^2$ and $10^2$, $T=3$ and $5$). In each replication, $k=500$ bootstrap samples are drawn to compute the p-value. In this simulation, our primary focus is on examining the test size and power of the proposed test. 

\subsubsection{Empirical sizes of the proposed test}

When $c=0$, $\beta_3(\tau_t)=-5$ and $\beta_4(\tau_t)=5$, meaning that the null hypothesis holds. For all experimental settings, given the significance levels \(\alpha = \{0.01, 0.05, 0.1\}\), the empirical sizes  presented in Table 3 are the  rejection rates of the proposed test  under the null hypothesis.
\begin{table}[h]
	\centering
	\footnotesize
	\vskip -0.3cm
	\caption{Empirical sizes of the proposed  test
 at significance levels $\alpha =0.01,0.05 $ and $0.10$}
	\label{table3}
        \renewcommand{\arraystretch}{1}  
	\resizebox{\textwidth}{!}{  
		\begin{tabular}{ccccccccc}\hline
			& & &\multicolumn{3}{c}{Rook}&\multicolumn{3}{c}{Queen}\\
			$\rho (\tau_t)$&Distribution of \(\varepsilon_{i t}\)&\((T,N)\)&
  $\alpha=0.01$&$\alpha=0.05$ &$\alpha=0.10$&$\alpha=0.01$&$\alpha=0.05$ &$\alpha=0.10$\\\hline
			 &\(N\left(0,1\right)\) &(3,64)&0.014&0.050&0.110&0.010&0.052&0.096\\
			&&(3,100)&0.018&0.074&0.140&0.010&0.068&0.124\\
			&&(5,64)&0.012&0.042&0.102&0.006&0.052&0.108\\
			&&(5,100)&0.004&0.042&0.088&0.024&0.052&0.092\\
		$\rho_1(\tau_t)$&$U\left(\sqrt{3},\sqrt{3}\right)$ & 
    (3,64)&0.008&0.046&0.088&0.004&0.046&0.098\\
			&&(3,100)&0.006&0.042&0.092&0.004&0.038&0.076\\
			&&(5,64)&0.010&0.056&0.096&0.012&0.046&0.090\\
			&&(5,100)&0.010&0.066&0.122&0.018&0.062&0.104\\
			&\(\frac{1}{2} \chi^{2}(2)-1\)
			&(3,64)&0.010&0.056&0.106&0.008&0.060&0.098\\
			&&(3,100)&0.006&0.060&0.086&0.004&0.046&0.102\\
			&&(5,64)&0.018&0.056&0.102&0.010&0.050&0.106\\
			&&(5,100)&0.012&0.050&0.110&0.026&0.066&0.110\\ 
    &\(N\left(0,1\right)\) &(3,64)&0.008&0.046&0.080&0.014&0.042&0.094\\
			&&(3,100)&0.012&0.070&0.122&0.008&0.076&0.132\\
			&&(5,64)&0.010&0.078&0.114&0.006&0.058&0.106\\
			&&(5,100)&0.010&0.074&0.134&0.010&0.036&0.098\\
		$\rho_2(\tau_t)$&$U\left(\sqrt{3},\sqrt{3}\right)$ &(3,64)&0.014&0.030&0.088&0.004&0.044&0.090\\
			&&(3,100)&0.010&0.034&0.092&0.008&0.032&0.074\\
			&&(5,64)&0.010&0.040&0.094&0.002&0.036&0.092\\
			&&(5,100)&0.020&0.072&0.132&0.006&0.044&0.082\\
			&\(\frac{1}{2} \chi^{2}(2)-1\)
			&(3,64)&0.008&0.046&0.080&0.014&0.042&0.094\\
			&&(3,100)&0.008&0.064&0.124&0.008&0.052&0.110\\
			&&(5,64)&0.016&0.054&0.094&0.008&0.036&0.076\\
			&&(5,100)&0.014&0.058&0.116&0.000&0.046&0.096\\\hline
		\end{tabular}
	}
\end{table}

As expected, all empirical sizes are close to the corresponding significance levels, indicating that the bootstrap method effectively approximates the null distribution of the test statistic. Furthermore, there is no significant difference in the rejection rates under the three different error term distributions, the two different spatial adjacency matrices and the two different forms of the spatial lag time-varying coefficient.  This suggests that the proposed test demonstrates robustness in maintaining the test size.

\subsubsection{Power of the test}
When $c \neq 0$, the alternative hypothesis that $\beta_3(\tau_t)$ and $\beta_4(\tau_t)$ are time-varying coefficients holds, the rejection rates are the estimates of the power of the test under the alternative hypothesis. Here, we set the the significance level $\alpha=0.05$ and the value of $c$ are drawn from the set $\{0.3,0.4,0.5\}$. We depict the power of the test with different value of $c$ under the three types of error distribution, the two form of the spatial lag coefficient and the two types of adjacent matrices in Figure 2.

\begin{figure}[htbp]
	\centering
	\begin{subfigure}{0.45\textwidth}
		\centering
		\includegraphics[width=\textwidth]{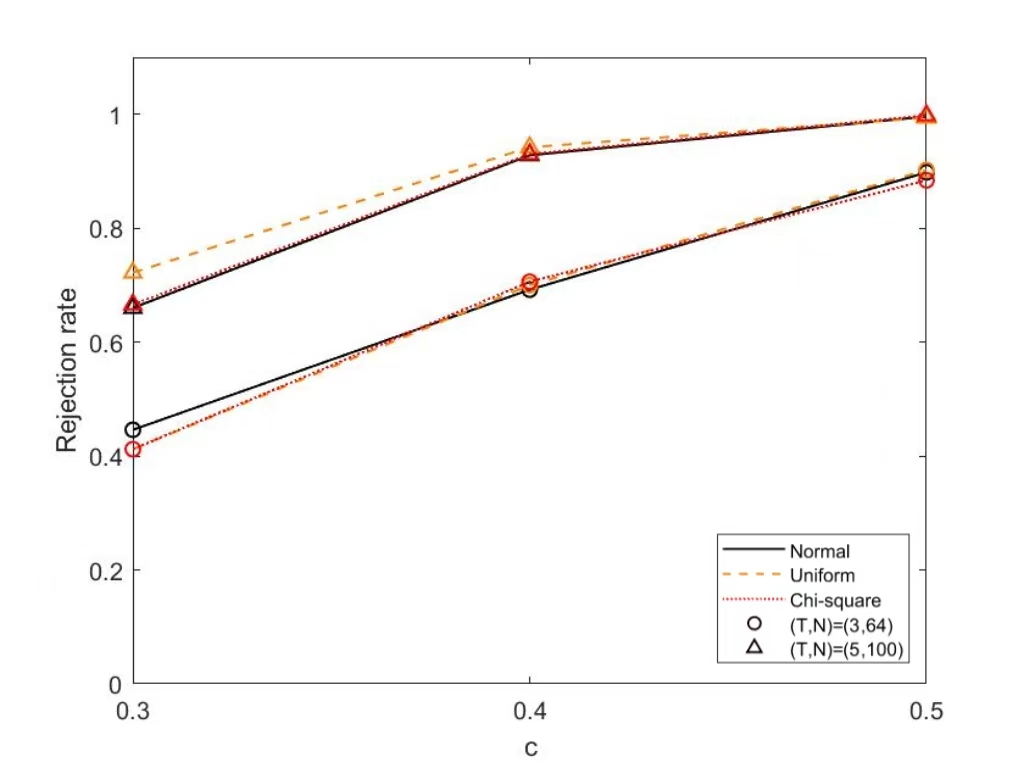}
		\caption{$\rho_1(\tau_t)$,Rook}
		\label{fig:subfig1}
	\end{subfigure}
	\hfill
	\begin{subfigure}{0.45\textwidth}
		\centering
		\includegraphics[width=\textwidth]{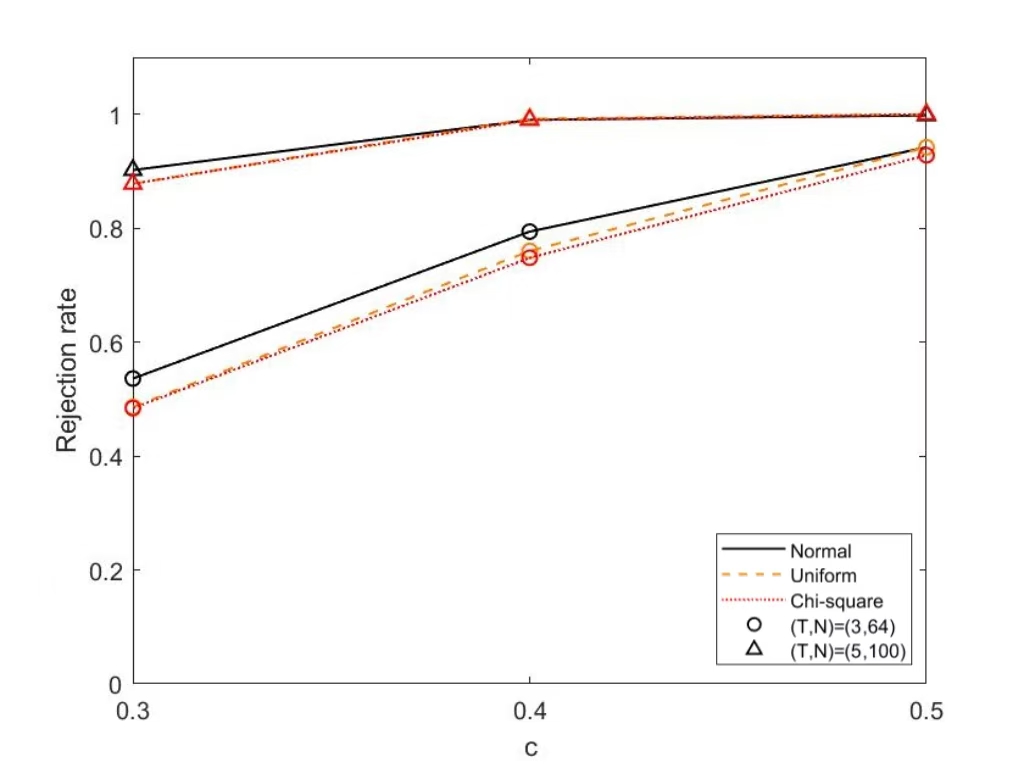}
		\caption{$\rho_1(\tau_t)$,Queen}
		\label{fig:subfig2}
	\end{subfigure}
	\begin{subfigure}{0.45\textwidth}
		\centering
		\includegraphics[width=\textwidth]{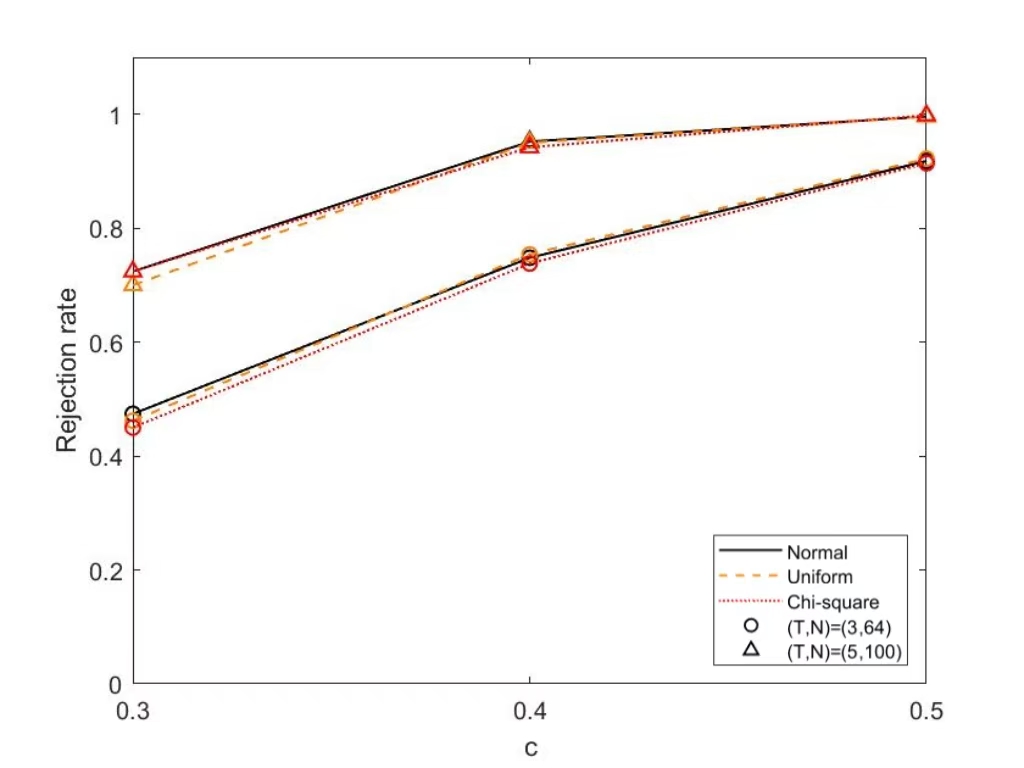}
		\caption{$\rho_2(\tau_t)$,Rook}
		\label{fig:subfig1}
	\end{subfigure}
	\hfill
	\begin{subfigure}{0.45\textwidth}
		\centering
		\includegraphics[width=\textwidth]{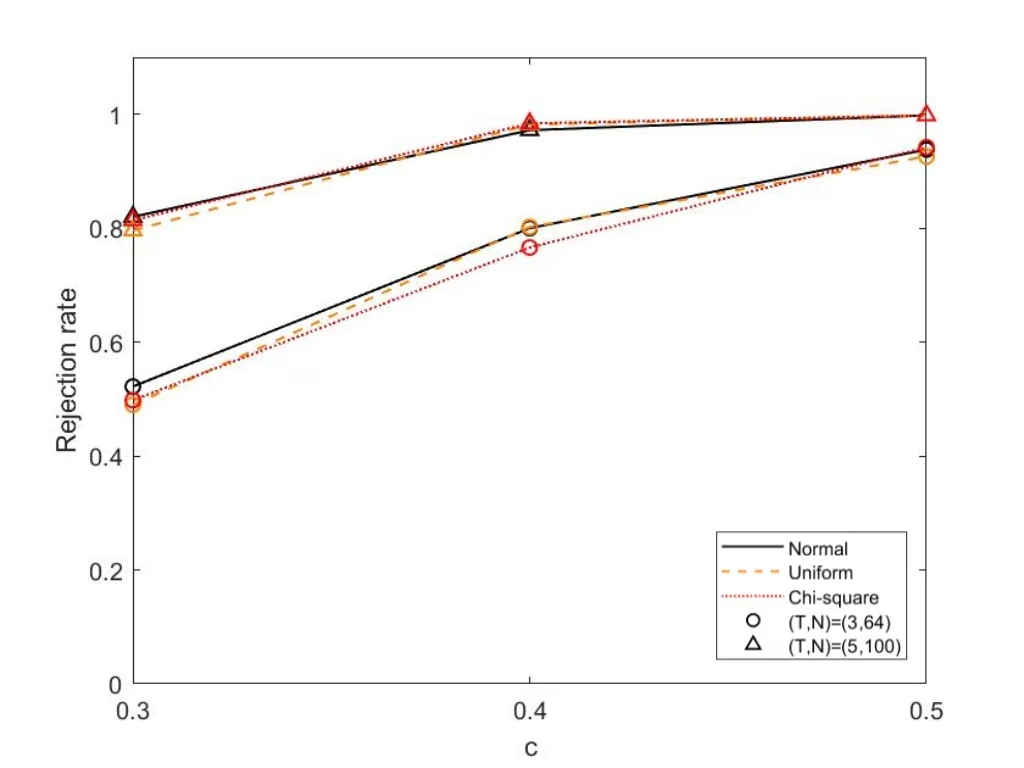}
		\caption{$\rho_2(\tau_t)$,Queen}
		\label{fig:subfig2}
	\end{subfigure}
	\caption{The power of the test under the significance level $\alpha =0.05$}
	\label{fig:twosubfigures}
\end{figure}

As shown in Figure 2, the test power rapidly approaches 1 across all experimental settings as the value of $c$ increases, indicating satisfactory performance. Notably, higher power is observed for larger values of $T$ or $N$. Additionally, the test maintains consistent power under various error distributions and different time-varying spatial autoregressive coefficients, demonstrating the robustness of the proposed method. Furthermore, the type of adjacency matrix significantly influences the rejection rates, with the test under the Queen adjacency matrix exhibiting significantly higher power compared to the Rook adjacency matrix, suggesting greater sensitivity under the Queen adjacency matrix.

\section{Empirical Application}
In this section, the proposed estimation and test are utilized to analyze provincial carbon emission data from China, showcasing their practical applications. The empirical dataset encompasses information from 30 provinces, autonomous regions, and municipalities directly administered by the Chinese central government (excluding Tibet, Hong Kong, Macao, and Taiwan), spanning the years 2005 to 2016. All panel data utilized in this study have been sourced from the National Bureau of Statistics' database.

For the purpose of comparison, we select the same four explanatory variables as those in Tian LL {\it et al.} (2024)  for our analysis. These   variables are respectively: \\

{\it PC} (Per capita carbon emissions): The ratio of the total carbon emissions and the the year-end resident population, with the unit being Tons per person. 

{\it PG} (Per capita GDP): The ratio of the total GDP and the resident population at the end of the year of each region, with the unit being Ten thousand yuan per person.

{\it PR} (Population structure): The ratio of urban population and total population.

{\it IR} (Industrial structure): The ratio of the added value of the secondary industry.

{\it ER} (Energy structure): The ratio of total energy consumption and regional GDP, with the unit being  Tons of standard coal per ten thousand yuan.\\

Tian  {\it et al.} (2024) built a TVC-SAR model between {\it PC} and the above four explanatory variables, where {\it PC} presents the per capita carbon emission, which can be calculate by the ratio of the total carbon emission and the resident population. They used the generalized likelihood ratio test based on the quasi-maximum likelihood estimation to detect the existence of the spatial autocorrelation and possible constant coefficients in the regression relationship, respectively. It was found that there is strong spatial autocorrelation in carbon emission and the coefficients of {\it PG} and {\it IR} are constant. Based on the results of Tian {\it et al.} (2024), we consider the partially linear time-varying coefficients spatial autoregressive panel data model with fixed effects:

\begin{equation} \label{20}
\begin{array}{rcl}
PC_{it} &=& \rho(\tau_t) \sum\limits_{j=1}^{30} w_{ij}PC_{jt} + \beta_0(\tau_t) + \beta_1(\tau_t)PR_{it} + \beta_2(\tau_t)ER_{it} \\
        & & + \beta_3PG_{it} + \beta_4IR_{it} + \alpha_i + \varepsilon_{it}, \quad i=1,2, \cdots, 30; \, t=1,2, \cdots, 12.
\end{array}
\end{equation}
where the Queen adjacency matrix is chosen as the spatial adjacency matrix, and the Gaussian kernel function with the ROT bandwidth were used to estimate the time-varying coefficients in estimation and testing process. 

\begin{figure}[htbp]
	\centering
	\begin{subfigure}{0.45\textwidth}
		\centering
		\includegraphics[width=\textwidth]{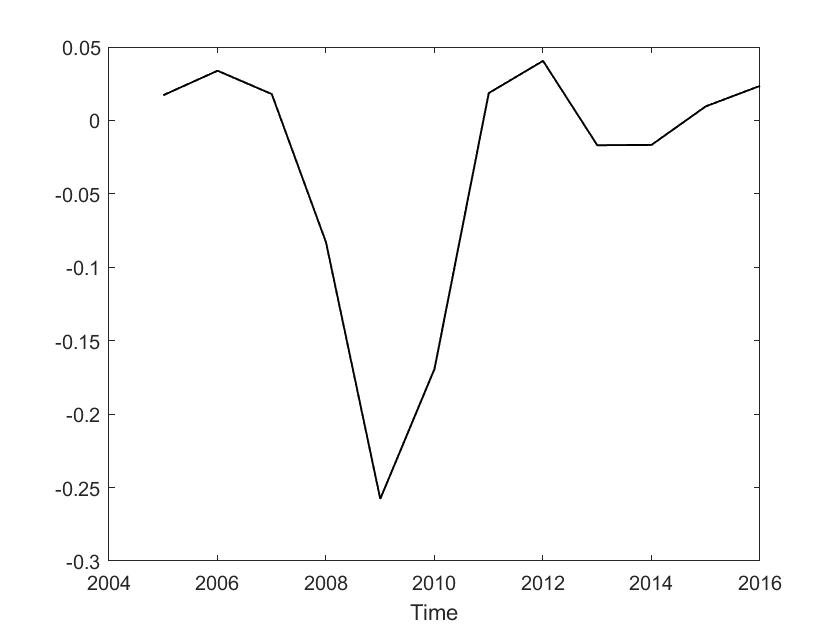}
		\caption{$\hat\rho(\tau_t)$}
		\label{fig:subfig1}
	\end{subfigure}
	\hfill
	\begin{subfigure}{0.45\textwidth}
		\centering
		\includegraphics[width=\textwidth]{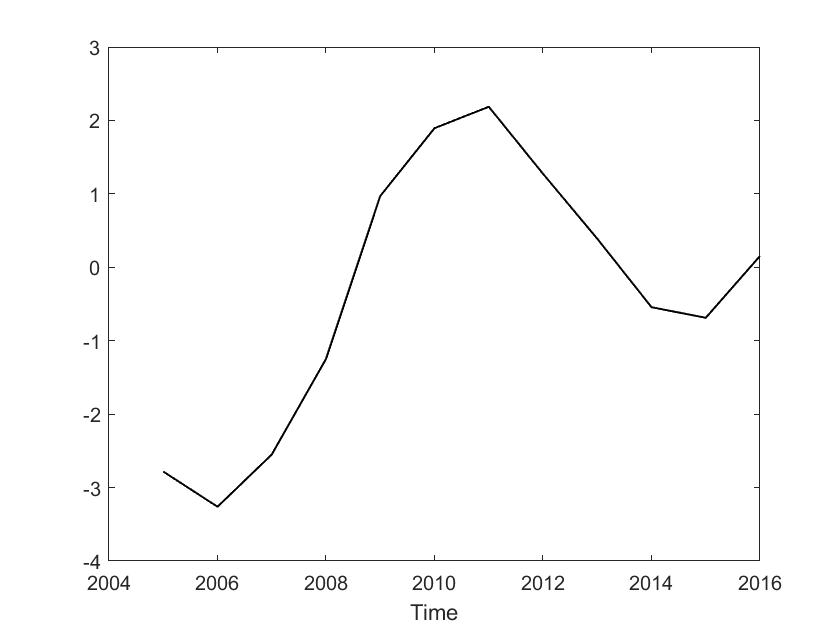}
		\caption{$\hat\beta_0(\tau_t)$}
		\label{fig:subfig2}
	\end{subfigure}
        \hfill
        \begin{subfigure}{0.45\textwidth}
		\centering
		\includegraphics[width=\textwidth]{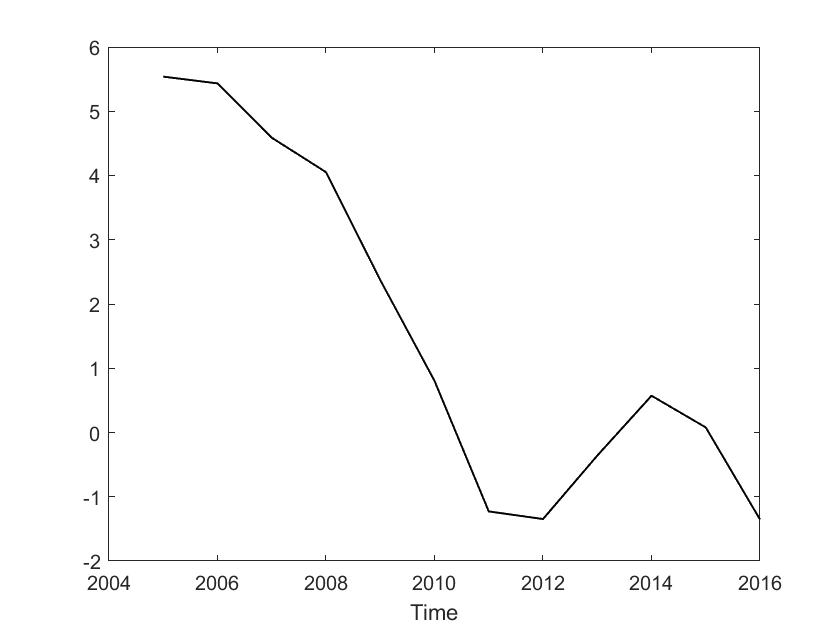}
		\caption{$\hat\beta_1(\tau_t)$}
		\label{fig:subfig3}
	\end{subfigure}
	\hfill
        \begin{subfigure}{0.45\textwidth}
		\centering
		\includegraphics[width=\textwidth]{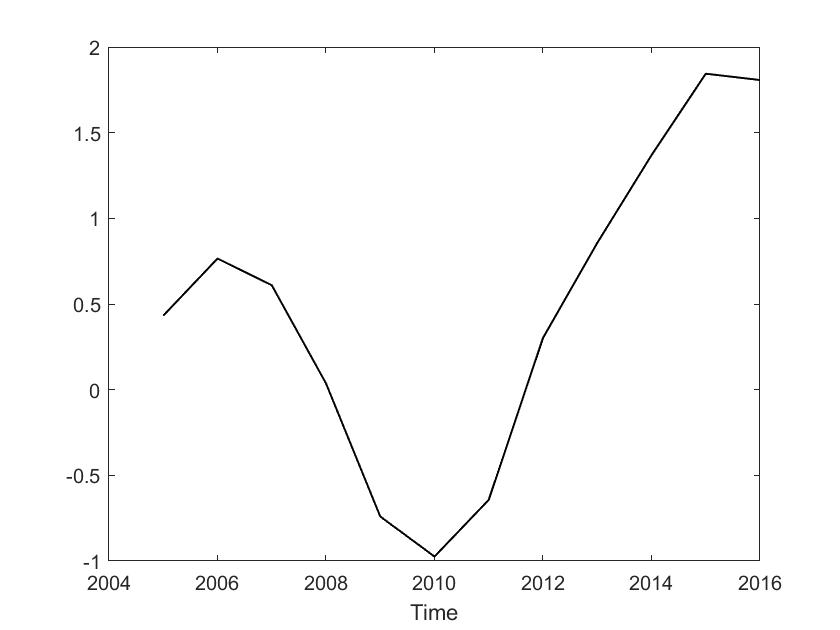}
		\caption{$\hat\beta_2(\tau_t)$}
		\label{fig:subfig4}
	\end{subfigure}
	\caption{The fitting curve of the time-varying coefficients in model (20).}
	\label{fig3}
\end{figure}

Firstly, we use the bootstrap test in Section 3 to check the hypothesis testing (3), a total of $k=500$ bootstrap samples were generated to calculate the p-value of the test. The resulting p-value is $0.204$, which provides further evidence that the coefficients of {\it PG} and {\it IR} are constant, Model (20) may be more appropriate for the carbon emission data. Then we employ our proposed 2SLS-PLLDV estimation in section 2 to fit Model (20), we obtain the estimates of the constant coefficients are $\hat\beta_3=0.1585$ and $\hat\beta_4=1.2111$, the fitting curves of the time-varying coefficients are illustrated in Figure \ref{fig3}.

In comparison to the regression coefficients reported by Tian  {\it et al.} (2024), where the coefficient of the spatial lag term remained constant, our proposed model identifies a broader range of time-varying coefficients. Notably, apart from the time-varying coefficients of {\it PR} (Population structure) and {\it ER} (Economic structure) as identified by Tian  {\it et al.} (2024), our model also pinpoints the spatial lag term and intercept coefficients as time-varying. The temporal  trend  of the spatial lag coefficient, depicted in Figure \ref{fig3}, highlights the dynamics of the autocorrelation coefficient. Initially, positive values signify spatial clustering effect,  indicating a positive correlation among high-carbon emission regions. Subsequently, possibly influenced by policy shifts, technological breakthroughs, or shifts in the economic landscape, the spatial autocorrelation coefficient transitions to negative, suggesting a compensatory effect among neighboring regions' carbon emissions. The subsequent resurgence of positive values could indicate the emergence of a fresh spatial clustering trend. Furthermore, the time-varying intercept coefficients offer a more nuanced understanding of how the baseline level of carbon emissions varies over time. Additionally, similar to Tian  {\it et al.} (2024), we find that the estimates for the constant coefficients of {\it PG}   and {\it IR} are positive, with consistent signs and comparable magnitudes, despite being derived from distinct methods.

\section{Conclusion}

We consider a partially linear time-varying coefficients SAR  panel data model that incorporates  fixed effects  and time-varying spatial lag coefficient. We have proposed the 2SLS-PLLDV method, which eliminates the need for first differencing to address fixed effects.   Under mild conditions, we establish the  asymptotic normality of the proposed estimators.  Additionally, we have developed  a residual-based test statistic to check the linear relationship between partial explanatory variables and the response variable, utilizing  a residual-based bootstrap procedure to calculate the p-values of the tests. Simulation studies  demonstrate the robust finite sample performance of both our proposed estimators and test statistic, exhibiting low  bias and high efficiency, as well as  reliable hypothesis testing outcomes. 
Additionally, we demonstrated the practical application of our method using provincial carbon emission data from China. Our empirical analysis revealed that the proposed method could capture more
complex time-varying relationships and offer insights into the temporal dynamics and spatial dependencies of carbon emissions.

The methods proposed in this paper not only enrich the statistical  method for SAR panel models but also possess substantial practical implications for model selection and policy analysis. By accommodating both constant and time-varying coefficients, our approach enhances the flexibility and robustness of spatial panel data modeling. Future endeavors could extend the model by exploring more adaptable functional forms to capture intricate nonlinear relationships. Additionally, comparative analyses with other advanced estimation techniques would further strengthen the validation of our method and identify areas for improvement.

\section*{CRediT authorship contribution statement}
Lingling Tian: Conceptualization, Methodology, Investigation, Formal analysis Writing-original draft.
Chuanhua Wei: Conceptualization, Methodology, Resources.
Mixia
Wu: Supervision, Methodology, Writing-review \& editing.

\section*{Acknowledgments}

The research is supported by  the Open Fund Project of Key Laboratory of Market Regulation (No: 2023SYSKF02003) and  the National Social Science Foundation of China (21BTJ005). 


\section*{Appendix. Proof of the main results}
Before providing the main results, we introduce the following useful lemmas.\\
{\bf \large Lemma 1.} Under    Assumptions 1-6, we have 
$$
\frac{1}{NTh} {\bf M}^{\top}(\tau_0){\bf W}_h(\tau_0){\bf M}(\tau_0)=\bm \Lambda_\mu \otimes\bm \Sigma_{{\bf z}_v}+o_p(1),
$$
where $\bm \Lambda_\mu=\text{diag}(\mu_0,\mu_2)$, and $\bm \Sigma_{{\bf z}_v}$ is as defined in Section 3.\\
{\bf Proof:} To prove this lemma, we first express the matrix ${\bf M}^\top(\tau_0) {\bf W}_h(\tau_0) {\bf M}(\tau_0)$  in its entirety. Notice that
$$
{\bf M}^{\top}(\tau_0){\bf W}_h(\tau_0){\bf M}(\tau_0)= \begin{pmatrix}
\sum\limits_{i=1}^N \sum\limits_{t=1}^T {\bf z}_{v,it} {\bf z}_{v,it}^\top K_t  & \sum\limits_{i=1}^N \sum\limits_{t=1}^T {\bf z}_{v,it} {\bf z}_{v,it}^\top \left(\frac{\tau_t-\tau_0}{h}\right) K_t\\
        \sum\limits_{i=1}^N \sum\limits_{t=1}^T {\bf z}_{v,it} {\bf z}_{v,it}^\top \left(\frac{\tau_t-\tau_0}{h}\right) K_t &\sum\limits_{i=1}^N \sum\limits_{t=1}^T {\bf z}_{v,it} {\bf z}_{v,it}^\top \left(\frac{\tau_t-\tau_0}{h} \right)^2 K_t
\end{pmatrix},
$$
where $K_t=K(\frac{\tau_t-\tau_0}{h})$. Since the proofs for the other components in ${\bf M}^{\top}(\tau_0){\bf W}_h(\tau_0){\bf M}(\tau_0)$ are analogous, we only need to prove 
$$
\frac{1}{NTh} \sum\limits_{i=1}^N \sum\limits_{t=1}^T {\bf z}_{v,it} {\bf z}_{v,it}^\top K_t=\mu_0\bm \Sigma_{{\bf z}_v}+o_p(1).     \eqno (\text{A.1})
$$

Firstly, we consider the expection of $\frac{1}{NTh} \sum\limits_{i=1}^N \sum\limits_{t=1}^T {\bf z}_{v,it} {\bf z}_{v,it}^\top K_t$. Based on  Assumptions 1-2, it is easy to show that 
$$
E\left(\frac{1}{NTh}{\sum\limits_{i=1}^N \sum\limits_{t=1}^T {\bf z}_{v,it} {\bf z}_{v,it}^\top K_t}\right) = \frac{1}{NTh} \sum\limits_{t=1}^T K_t \sum\limits_{i=1}^N E\left({\bf z}_{v,it} {\bf z}_{v,it}^\top\right).   
$$
Since $E({\bf z}_{v,it} {\bf z}_{v,it}^\top) = \bm\Sigma_{{\bf z}_v}$  and $\frac{1}{Th} \sum_{t=1}^T K_t \to \mu_0$  as both \(T \to \infty\) and \(Th \to \infty\), we have
$$
E\left(\frac{1}{NTh}{\sum\limits_{i=1}^N \sum\limits_{t=1}^T {\bf z}_{v,it} {\bf z}_{v,it}^\top K_t}\right) = \mu_0 \bm\Sigma_{{\bf z}_v}\left(1+O(\frac{1}{Th})\right).  \eqno (\text{A.2})
$$

Then we consider the variance of $\frac{1}{NTh} \sum\limits_{i=1}^N \sum\limits_{t=1}^T {\bf z}_{v,it} {\bf z}_{v,it}^\top K_t$. By following the proofs of (A.5)-(A.12) in Li {\it et al.} (2011), we can deduce
$$
{\rm Var}\left(\frac{1}{NTh}{\sum\limits_{i=1}^N \sum\limits_{t=1}^T {\bf z}_{v,it} {\bf z}_{v,it}^\top K_t}\right)=O\left(\frac{1}{NTh}\right).  \eqno (\text{A.3})
$$
By combining the results (A.2) and (A.3), we have shown that (A.1) holds. The proof of Lemma 1 is completed.\\ 
 {\bf \large Lemma 2.} Under Assumptions 1-6, we have
 $$
 \frac{1}{NT}{\bar{\bf X}}_c^\top {\bar{\bf X}}_c \overset{P}{\rightarrow} \bm\Sigma,
 $$
 where $\bm\Sigma=\bm\Sigma_{{\bf x}_c}-\bm\Sigma_{{\bf z}_v {\bf x}_c}^\top \bm\Sigma_{{\bf z}_v}^{-1}\bm\Sigma_{{\bf z}_v {\bf x}_c}$.\\
 {\bf Proof:} Upon expanding the expression of ${\bar{\bf X}}_c^\top {\bar{\bf X}}_c$, we obtain
 $$
 \begin{array}{ll}
      & {\bar{\bf X}}_c^\top {\bar{\bf X}}_c\\
     =&{\bf X}_c^\top ({\bf I}_{NT} - {\bf S}(h))^\top ({\bf I}_{NT} - {\bf P}_D) ({\bf I}_{NT} - {\bf S}(h)) {\bf X}_c \\
     = & {\bf X}_c^\top ({\bf I}_{NT} - {\bf S}(h))^\top ({\bf I}_{NT} - {\bf S}(h)) {\bf X}_c-{\bf X}_c^\top ({\bf I}_{NT} - {\bf S}(h))^\top {\bf D}({\bf D}^\top {\bf D})^{-1}{\bf D}^\top ({\bf I}_{NT} - {\bf S}(h)) {\bf X}_c\\
     =&\mathbb{A}_1 + \mathbb{A}_2.   
 \end{array}     \eqno (\text{A.4})
 $$
We analyze each term individually. According to Lemma 3 and Lemma 4 in   in Zhao {\it et al.} (2016), we have $\mathbb{A}_2=o_P(1)$.

For $\mathbb{A}_1$, we find that
 $$
 \begin{array}{ll}
      & ({\bf I}_{NT} - {\bf S}(h)) {\bf X}_c \\
     = & \sum \limits_{i=1}^N \sum\limits_{t=1}^T \left\{{\bf x}^\top_{c,it}-({\bf z}_{v,it}^\top,{\bf 0}_{1\times(q+1)})\left[{\bf M}^{\top}(\tau_0){\bf W}_h^{*}(\tau_0){\bf M}(\tau_0)\right]^{-1} {\bf M}^{\top}(\tau_0){\bf W}_h^{*}(\tau_0){\bf X}_c\right\}. 
 \end{array}
 $$
By the definition of ${\bf W}_h^{*}(\tau_0)$, we have
$$
\begin{array}{ll}
     & {\bf M}^{\top}(\tau_0){\bf W}_h^{*}(\tau_0){\bf M}(\tau_0) \\
     =& {\bf M}^{\top}(\tau_0){\bf W}_h(\tau_0){\bf M}(\tau_0)-{\bf M}^{\top}(\tau_0){\bf W}_h(\tau_0){\bf D}\left({\bf D}^\top {\bf W}_h(\tau_0) {\bf D}\right)^{-1}{\bf D}^\top {\bf W}_h(\tau_0){\bf M}(\tau_0),  
\end{array}  \eqno (\text{A.5})
$$
and 
$$
\begin{array}{ll}
     & {\bf M}^{\top}(\tau_0){\bf W}_h^{*}(\tau_0){\bf X}_c \\
     =& {\bf M}^{\top}(\tau_0){\bf W}_h(\tau_0){\bf X}_c -{\bf M}^{\top}(\tau_0){\bf W}_h(\tau_0){\bf D}\left({\bf D}^\top {\bf W}_h(\tau_0) {\bf D}\right)^{-1}{\bf D}^\top {\bf W}_h(\tau_0){\bf X}_c.   
\end{array}   \eqno (\text{A.6})
$$
By following the proof of (A.26) in Li {\it et al.} (2011), we obtain 
$$
{\bf M}^{\top}(\tau_0){\bf W}_h(\tau_0){\bf D}\left({\bf D}^\top {\bf W}_h(\tau_0) {\bf D}\right)^{-1}{\bf D}^\top {\bf W}_h(\tau_0){\bf M}(\tau_0)=o_P(NTh),  \eqno (\text{A.7})
$$
$$
{\bf M}^{\top}(\tau_0){\bf W}_h(\tau_0){\bf D}\left({\bf D}^\top {\bf W}_h(\tau_0) {\bf D}\right)^{-1}{\bf D}^\top {\bf W}_h(\tau_0){\bf X}_c =o_P(NTh).  \eqno (\text{A.8})
$$
Together with  Lemma 1, we conclude that $\frac{1}{NTh} {\bf M}^{\top}(\tau_0){\bf W}_h^{*}(\tau_0){\bf M}(\tau_0)=\bm\Lambda_\mu \otimes\bm\Sigma_{{\bf z}_v}+o_P(1)$. Now we need only to prove 
$$
\frac{1}{NTh}{\bf M}^{\top}(\tau_0){\bf W}_h(\tau_0){\bf X}_c=\begin{pmatrix}
    1\\
    0
\end{pmatrix} \otimes\bm\Sigma_{{\bf z}_v{\bf x}_c}+o_P(1).   \eqno (\text{A.9}) 
$$
Similar to the proof of (A.1) in Lemma 1, we have 
$$
E\left[\frac{1}{NTh}{\bf M}^{\top}(\tau_0){\bf W}_h(\tau_0){\bf X}_c\right]=\frac{1}{NTh}\begin{pmatrix}
    \sum\limits_{i=1}^N  \sum\limits_{t=1}^T K_tE({\bf z}_{v,it}{\bf x}^\top_{c,it})\\
    \sum\limits_{i=1}^N  \sum\limits_{t=1}^T (\frac{\tau_t-\tau_0}{h})K_tE({\bf z}_{v,it}{\bf x}^\top_{c,it})
\end{pmatrix}=\begin{pmatrix}
        \bm\Sigma_{{\bf z}_v{\bf x}_c}\\
        0
    \end{pmatrix},
$$
$$
Var\left[\frac{1}{NTh}{\bf M}^{\top}(\tau_0){\bf W}_h(\tau_0){\bf X}_c\right]=O(\frac{1}{NTh}).
$$
Thus,  (A.9) holds. Through straightforward calculation, we derive
$$
({\bf z}_{v,it}^\top,{\bf 0}_{1\times(q+1)})\left[{\bf M}^{\top}(\tau_0){\bf W}_h^{*}(\tau_0){\bf M}(\tau_0)\right]^{-1} {\bf M}^{\top}(\tau_0){\bf W}_h^{*}(\tau_0){\bf X}_c={\bf z}_{v,it}^\top\bm\Sigma_{{\bf z}_v}^{-1}\Sigma_{{\bf z}_v{\bf x}_c},
$$
and 
$$
\begin{array}{ll}
     &E\left\{\frac{1}{NT}{\bf X}_c^\top ({\bf I}_{NT} - {\bf S}(h))^\top ({\bf I}_{NT} - {\bf S}(h)) {\bf X}_c\right\}\\
     =&\frac{1}{NT}\sum \limits_{i=1}^N \sum\limits_{t=1}^T E\left\{{\bf x}^\top_{c,it}-{\bf z}_{v,it}^\top\bm\Sigma_{{\bf z}_v}^{-1}\bm\Sigma_{{\bf z}_v{\bf x}_c}\right\}^\top \left\{{\bf x}^\top_{c,it}-{\bf z}_{v,it}^\top\bm\Sigma_{{\bf z}_v}^{-1}\bm\Sigma_{{\bf z}_v{\bf x}_c}\right\}\\
     =&\bm\Sigma_{{\bf x}_c}-\bm\Sigma^\top_{{\bf z}_v{\bf x}_c}\bm\Sigma_{{\bf z}_v}^{-1}\bm\Sigma_{{\bf z}_v{\bf x}_c}. 
\end{array}
$$
According to the Slutsky Theorem, we have  
$$
\frac{1}{NT}{\bar{\bf X}}_c^\top {\bar{\bf X}}_c \overset{P}{\rightarrow}\bm\Sigma_{{\bf x}_c}-\bm\Sigma^\top_{{\bf z}_v{\bf x}_c}\bm\Sigma_{{\bf z}_v}^{-1}\bm\Sigma_{{\bf z}_v{\bf x}_c}.
$$
Hence, Lemma 2 holds.\\
{{\bf \large Lemma 3.} } Under  Assumptions 1-6, we have
 $$
 \frac{1}{\sqrt{NT}}{\bf X}_c^\top ({\bf I}_{NT} - {\bf S}(h))^\top ({\bf I}_{NT} - {\bf S}(h)) {\bm \varepsilon} \overset{D}{\rightarrow} N({\bf 0}_{p-q},\bm\Omega_{\bm \varepsilon})
 $$
 where $\bm\Omega_{\bm\varepsilon}=\sum\limits_{t=-\infty}^\infty c_{\bm\varepsilon}(t)\left({\bf C}_{{\bf x}_c}(t)-{\bf C}_{{{\bf x}_c}{{\bf z}_v}}(t)\bm\Sigma_{{\bf z}_v}^{-1} \bm\Sigma_{{\bf z}_v {\bf x}_c}  +\bm\Sigma^\top_{{\bf z}_v {\bf x}_c}\bm\Sigma_{{\bf z}_v}^{-1}{\bf C}_{{\bf z}_v}(t)\bm\Sigma_{{\bf z}_v}^{-1} \bm\Sigma_{{\bf z}_v {\bf x}_c}   \right.$\\
 $\left.-\bm\Sigma^\top_{{\bf z}_v {\bf x}_c}\bm\Sigma_{{\bf z}_v}^{-1}{\bf C}_{{{\bf z}_v}{{\bf x}_c}}(t)\right). 
 $ 
\\

 {\bf Proof:} Note that
 $$
 \begin{array}{ll}
      &  {\bf X}_c^\top ({\bf I}_{NT} - {\bf S}(h))^\top ({\bf I}_{NT} - {\bf S}(h)) {\bm \varepsilon} \\
     = & \sum \limits_{i=1}^N \sum\limits_{t=1}^T \left\{{\bf x}^\top_{c,it}-{\bf z}_{v,it}^\top\bm\Sigma_{{\bf z}_v}^{-1}\bm\Sigma_{{{\bf z}_v}{{\bf x}_c}}\right\}^\top \left\{{\varepsilon}_{it}-{\bf z}_{v,it}^\top\bm\Sigma_{{\bf z}_v}^{-1}{\bf M}^{\top}(\tau_0){\bf W}_h^{*}(\tau_0){\bm \varepsilon}\right\}.
 \end{array}  \eqno (\text{A.10})
 $$
 Analogous to the derivations in (A.5)-(A.8), we have ${\bf M}^{\top}(\tau_0){\bf W}_h^{*}(\tau_0){\bm \varepsilon}={\bf M}^{\top}(\tau_0){\bf W}_h(\tau_0){\bm \varepsilon}+o_P(NTh)$. Since
 $$
 E\left(\frac{1}{NTh}{\bf M}^{\top}(\tau_0){\bf W}_h(\tau_0){\bm \varepsilon}\right)=\frac{1}{NTh}E\begin{pmatrix}
       \sum\limits_{i=1}^N  \sum\limits_{t=1}^T K_t{\bf z}_{v,it}{\varepsilon}_{it}\\
    \sum\limits_{i=1}^N  \sum\limits_{t=1}^T (\frac{\tau_t-\tau_0}{h})K_t{\bf z}_{v,it}{\varepsilon}_{it}
 \end{pmatrix}=o(1),
 $$
we can simplify (A.10)  to
 $$
 {\bf X}_c^\top ({\bf I}_{NT} - {\bf S}(h))^\top ({\bf I}_{NT} - {\bf S}(h)) {\bm \varepsilon} = \sum \limits_{i=1}^N \sum\limits_{t=1}^T \left\{{\bf x}^\top_{c,it}-{\bf z}^\top_{v,it}\bm\Sigma_{{\bf z}_v}^{-1}\bm\Sigma_{{\bf z}_v{\bf x}_c}\right\}^\top{\varepsilon}_{it}+o_P(1).  \eqno (\text{A.11})
 $$
 The expectation and covariance matrix of (A.11) are as follows
 $$
 E\left[{\bf X}_c^\top ({\bf I}_{NT} - {\bf S}(h))^\top ({\bf I}_{NT} - {\bf S}(h)) {\bm \varepsilon}\right]=  o(1),
 $$
 $$
 \begin{array}{ll}
      &\frac{1}{NT}{\rm Var}\left[{\bf X}_c^\top ({\bf I}_{NT} - {\bf S}(h))^\top ({\bf I}_{NT} - {\bf S}(h)) {\bm \varepsilon}\right] \\
    =&\frac{1}{NT} \sum\limits_{i=1}^N \sum\limits_{j=1}^N \sum\limits_{t=1}^T \sum\limits_{s=1}^T E\left\{({\bf x}^\top_{c,it}-{\bf z}_{v,it}^\top\bm\Sigma_{{\bf z}_v}^{-1}\bm\Sigma_{{\bf z}_v{\bf x}_c})^\top({\bf x}^\top_{c,js}-{\bf z}_{v,js}^\top\bm\Sigma_{{\bf z}_v}^{-1}\bm\Sigma_{{\bf z}_v{\bf x}_c})\varepsilon_{it}\varepsilon_{js}\right\}\\
    =&\frac{1}{NT}\sum\limits_{i=1}^N \sum\limits_{t=1}^T \sum\limits_{s=1}^T E\left\{({\bf x}^\top_{c,it}-{\bf z}_{v,it}^\top\bm\Sigma_{{\bf z}_v}^{-1}\bm\Sigma_{{\bf z}_v{\bf x}_c})^\top({\bf x}^\top_{c,is}-{\bf z}_{v,is}^\top\bm\Sigma_{{\bf z}_v}^{-1}\bm\Sigma_{{\bf z}_v{\bf x}_c})\varepsilon_{it}\varepsilon_{is}\right\}\\
    =&\frac{1}{NT}\sum\limits_{i=1}^N \sum\limits_{t=1}^T\sum\limits_{s=t}\sigma_{\bm\varepsilon}^2\left\{\bm\Sigma_{{\bf x}_c}-\bm\Sigma_{{\bf z}_v {\bf x}_c}^\top \bm\Sigma_{{\bf z}_v}^{-1}\bm\Sigma_{{\bf z}_v {\bf x}_c}\right\} 
    +\sum\limits_{i=1}^N \sum\limits_{t=1}^T\sum\limits_{s\neq t}  c_{\bm\varepsilon}(s-t)\left\{{\bf C}_{{\bf x}_c}(s-t)\right.\\&
    \left.-{\bf C}_{{{\bf x}_c}{{\bf z}_v}}(s-t)\bm\Sigma_{{\bf z}_v}^{-1} \bm\Sigma_{{\bf z}_v {\bf x}_c}-\bm\Sigma^\top_{{\bf z}_v {\bf x}_c}\bm\Sigma_{{\bf z}_v}^{-1}{\bf C}_{{{\bf z}_v}{{\bf x}_c}}(s-t)+\bm\Sigma^\top_{{\bf z}_v {\bf x}_c}\bm\Sigma_{{\bf z}_v}^{-1}{\bf C}_{{\bf z}_v}(s-t)\bm\Sigma_{{\bf z}_v}^{-1} \bm\Sigma_{{\bf z}_v {\bf x}_c}\right\}\\
    =&\sum\limits_{t=1-T}^{T-1} c_{\bm\varepsilon}(t)\left\{{\bf C}_{{\bf x}_c}(t)-{\bf C}_{{{\bf x}_c}{{\bf z}_v}}(t)\bm\Sigma_{{\bf z}_v}^{-1} \bm\Sigma_{{\bf z}_v {\bf x}_c}-\bm\Sigma^\top_{{\bf z}_v {\bf x}_c}\bm\Sigma_{{\bf z}_v}^{-1}{\bf C}_{{{\bf z}_v}{{\bf x}_c}}(t)\right.\\
   &\left.+\bm\Sigma^\top_{{\bf z}_v {\bf x}_c}\bm\Sigma_{{\bf z}_v}^{-1}{\bf C}_{{\bf z}_v}(t)\bm\Sigma_{{\bf z}_v}^{-1} \bm\Sigma_{{\bf z}_v {\bf x}_c}\right\}\\
   \rightarrow & \bm\Omega_{\bm\varepsilon}\ 
   \mbox{\  (as $T\rightarrow \infty$).} 
   \end{array}$$ 
  Thus, Lemma 3 holds.\\
\textbf{Proof of Theorem 3.1.} We start by examining the expression for \(\hat{\bm{\beta}}_c - \bm{\beta}_c\):
\begin{align*}
    &\hat{\bm \beta}_c-{\bm \beta}_c=(\bar{\bf X}_c^\top \bar{\bf X}_c)^{-1}\bar{\bf X}_c^\top \bar {\bf Y}-{\bm \beta}_c\\
    &=(\bar{\bf X}_c^\top \bar{\bf X}_c)^{-1}\bar{\bf X}_c^\top({\bf I}_{NT}-{\bf P}_{\bf D})({\bf I}_{NT}-{\bf S}(h))\left[{\bf B}({\bf Z}_v,{\bm \gamma}_v)+{\bf X}_c\boldsymbol\beta_c+{\bf D} \boldsymbol\alpha+\boldsymbol\varepsilon\right]-{\bm \beta}_c \\
    &=(\bar{\bf X}_c^\top \bar{\bf X}_c)^{-1}\bar{\bf X}_c^\top({\bf I}_{NT}-{\bf P}_{\bf D})({\bf I}_{NT}-{\bf S}(h)){\bf B}({\bf Z}_v,{\bm \gamma}_v) \\
    &\quad +(\bar{\bf X}_c^\top \bar{\bf X}_c)^{-1}\bar{\bf X}_c^\top({\bf I}_{NT}-{\bf P}_{\bf D})({\bf I}_{NT}-{\bf S}(h)){\bf D} \boldsymbol\alpha \\
    &\quad +(\bar{\bf X}_c^\top \bar{\bf X}_c)^{-1}\bar{\bf X}_c^\top({\bf I}_{NT}-{\bf P}_{\bf D})({\bf I}_{NT}-{\bf S}(h))\boldsymbol\varepsilon \\
    &= \mathbb{D}_1 + \mathbb{D}_2 + \mathbb{D}_3.
\end{align*}
Based on the Lemma 4.5 from Su and Ullah (2006), we have 
$$
\frac{1}{\sqrt{NT}}({\bf I}_{NT}-{\bf P}_{\bf D})({\bf I}_{NT}-{\bf S}(h)){\bf B}({\bf Z}_v,{\bm \gamma}_v)=o_p(1),
$$
which implies that   $\mathbb{D}_1=o_p(1)$. For $\mathbb{D}_2$, by the definition of ${\bf W}_h^{*}(\tau_0)$ and ${\bf P}_{\bf D}$, we have ${\bf S}(h){\bf D}={\bf 0}_{NT\times(N-1)}$, $({\bf I}_{NT}-{\bf P}_{\bf D}){\bf D}={\bf 0}_{NT\times(N-1)}$.
Through straightforward calculation, we obtain   $({\bf I}_{NT}-{\bf P}_{\bf D})({\bf I}_{NT}-{\bf S}(h)){\bf D}={\bf 0}_{NT\times(N-1)}$, indicating that$\mathbb{D}_2={\bf 0}_{p-q}$. 
For $\mathbb{D}_3$, by performing  the same decomposition as in (A.4), we have
\begin{align*}
&(\bar{\bf X}_c^\top \bar{\bf X}_c)^{-1}\bar{\bf X}_c^\top({\bf I}_{NT}-{\bf P}_{\bf D})({\bf I}_{NT}-{\bf S}(h))\boldsymbol\varepsilon\\
=&(\bar{\bf X}_c^\top \bar{\bf X}_c)^{-1}{\bf X}_c^\top({\bf I}_{NT}-{\bf S}(h))^\top({\bf I}_{NT}-{\bf S}(h))\boldsymbol\varepsilon+o_P(1).
\end{align*}
According to Lemma 2 and Lemma 3, we have
$$
\sqrt{NT}\mathbb{D}_3 \overset{d}{\rightarrow} N({\bf 0}_{p-q},{\bm\Sigma}^{-1}{\bm\Omega}_{\bm\varepsilon}{\bm\Sigma}^{-1}).
$$\\
By combining the results of \(\mathbb{D}_1\), \(\mathbb{D}_2\), and \(\mathbb{D}_3\), we obtain
$$
\sqrt{NT}(\hat{\bm{\beta}}_c - \bm{\beta}_c) \overset{D}{\rightarrow} N({\bf 0}_{p-q}, {\bm\Sigma}^{-1}{\bm\Omega}_{\bm\varepsilon}{\bm\Sigma}^{-1}).
$$
Therefore, the estimator \(\hat{\bm{\beta}}_c\) is asymptotically normal, with a mean \(\bm{\beta}_c\) and a variance-covariance matrix given by ${\bm\Sigma}^{-1}{\bm\Omega}_{\bm\varepsilon}{\bm\Sigma}^{-1}$. This completes the proof of the asymptotic properties for \(\bm{\beta}_c\).\\
\textbf{Proof of Theorem 3.2.} Note that
\begin{align*}
    &\hat{\boldsymbol\gamma}_v(\tau_0)-{\boldsymbol\gamma}_v(\tau_0)={\bm{\Phi}}(\tau_0)({\bf Y}-{\bf X}_c \hat{ \boldsymbol\beta}_c)-{\boldsymbol\gamma}_v(\tau_0) \\
    &= {\bm{\Phi}}(\tau_0)\left[{\bf B}({\bf Z}_v,{\bm \gamma}_v)+{\bf X}_c(\boldsymbol\beta_c-\hat{ \boldsymbol\beta}_c)+{\bf D} \boldsymbol\alpha+\boldsymbol\varepsilon\right]-{\boldsymbol\gamma}_v(\tau_0) \\
    &= \left\{\left({\bf I}_{q+1}, {\bf 0 }_{(q+1) \times (q+1)}\right) \left\{{\bf M}^{\top}(\tau_0){\bf W}_h^{*}(\tau_0){\bf M}(\tau_0)\right\}^{-1} {\bf M}^{\top}(\tau_0){\bf W}_h^{*}(\tau_0){\bf B}({\bf Z}_v,{\bm \gamma}_v)-{\boldsymbol\gamma}_v(\tau_0)\right\} \\
    &\quad +\left({\bf I}_{q+1}, {\bf 0 }_{(q+1) \times (q+1)}\right) \left\{{\bf M}^{\top}(\tau_0){\bf W}_h^{*}(\tau_0){\bf M}(\tau_0)\right\}^{-1} {\bf M}^{\top}(\tau_0){\bf W}_h^{*}(\tau_0){\bf X}_c(\hat{ \boldsymbol\beta}_c-{ \boldsymbol\beta}_c) \\
    &\quad +\left({\bf I}_{q+1}, {\bf 0 }_{(q+1) \times (q+1)}\right) \left\{{\bf M}^{\top}(\tau_0){\bf W}_h^{*}(\tau_0){\bf M}(\tau_0)\right\}^{-1} {\bf M}^{\top}(\tau_0){\bf W}_h^{*}(\tau_0){\bf D} \boldsymbol\alpha \\
    &\quad +\left({\bf I}_{q+1}, {\bf 0 }_{(q+1) \times (q+1)}\right) \left\{{\bf M}^{\top}(\tau_0){\bf W}_h^{*}(\tau_0){\bf M}(\tau_0)\right\}^{-1} {\bf M}^{\top}(\tau_0){\bf W}_h^{*}(\tau_0)\boldsymbol\varepsilon \\
    &= \mathbb{E}_1 + \mathbb{E}_2 + \mathbb{E}_3 + \mathbb{E}_4.
\end{align*}

By Assumption 4 and the Taylor expansion of ${\boldsymbol\gamma}_v(\tau_t)$ at $\tau_0$, we have
$$
{\bf B}({\bf Z}_v,{\bm \gamma}_v) = {\bf M}(\tau_0)\left({\boldsymbol\gamma}^\top_v(\tau_0), h{{\boldsymbol\gamma}_v'}^\top(\tau_0)\right)^\top + \frac{1}{2} \mu_2 {\bm\gamma}_v''(\tau_0)h^2 + o_p(h^2),
$$
thus $\mathbb{E}_1 = \frac{1}{2} \mu_2 {\bm\gamma}_v''(\tau_0)h^2 + o_p(h^2)$.

It is easy to show that $\sqrt{NT}\mathbb{E}_2 = o_p(1)$ by using the results in Theorem 3.1. For $\mathbb{E}_3$, by the definition of ${\bf W}_h^{*}(\tau_0)$, we have ${\mathcal{K}}(\tau_0){\bf D}={\bf 0}_{NT\times(N-1)}$, hence $\mathbb{E}_3 = {\bf 0}_{q+1}$.

Next we foucs on $\mathbb{E}_4$. By combining with (A.5), (A.7) and Lemma 1, we derive the following equation:
$$
\frac{1}{NTh}{\bf M}^{\top}(\tau_0){\bf W}_h^{*}(\tau_0){\bf M}(\tau_0)=\Lambda_\mu \otimes {\bf\Sigma}_{{\bf Z}_v}+o_p(1),  \eqno (\text{A.12})
$$
Following Proposition (A.6) in Li {\it et al.} (2011), we obtain
 $$
\frac{1}{\sqrt{NTh}}{\bf M}^{\top}(\tau_0){\bf W}_h^{*}(\tau_0){\bm\varepsilon}\xrightarrow{d} N(\mathbf{0}_{2(q+1)}, \Lambda_\nu \otimes {\bf\Omega}_v),  \eqno (\text{A.13})
$$
where $\Lambda_\nu = \mathrm{diag}(\nu_0, \nu_2)$, ${\bf\Omega}_v =\sum\limits_{t=-\infty}^{\infty}{\bf C}_{{\bf z}_v}(t)c_{\bm \varepsilon}(t)$.
Using (A.12) and (A.13), we deduce that:
$$
\sqrt{NTh}\mathbb{E}_4 \xrightarrow{D} N(\mathbf{0}_{q+1}, {\bm\Sigma}_{{\bf z}_v}^{-1} {\bf \Omega} {\bf \Sigma}_{{\bf z}_v}^{-1}).
$$
where ${\bf \Omega}=\nu_0 {\bf\Omega}_v$. Based on the above discussion and  utilizing the Slutsky Theorem, we have proven  Theorem 3.2.

\section*{MATLAB Code}
\subsection*{Estimation}
\begin{lstlisting}
function Estimation

diary './Estimation:error=1T=5m=10c=queen.txt';
disp("program starts");

% Parameter Initialization
T=10;
m=12; 
N=m^2;
nsim=500;
rand('seed', 0);
randn('seed', 1);

% Generating row-standardized spatial weight matrix
w=generationrook(m,N);
%w=generationqueen(m,N);

%row_standardization
sumc=sum(w,2); %sum by row
dd=ones(N,1)./sumc;
d=diag(dd);
ww=d*w;   %row normalization
W=kron(eye(T),ww);    

% Begin Simulation
for i=1:nsim

    taut=(1:T)'./T;
    tau=kron(taut,ones(N,1));
    
    x1=ones(N*T,1);
    x2=randn(N*T,1);
    x3=randn(N*T,1);
    x4=randn(N*T,1);
    X=[x1 x2 x3 x4];
    X_v=[x1 x2];
    X_c=[x3 x4];
    
    rhot=-0.6*(sin(2*pi*taut).^2);
    %rhot=0.6*(sin(2*pi*taut).^2);
    rho_NT=kron(diag(rhot),eye(N));
    beta1t=4*taut;
    beta1=4*tau;
    beta2t=(taut+1).^2;
    beta2=(tau+1).^2;
    beta3=-5;
    beta4=5;

    alpha=rand(N-1,1); 
    D=kron(ones(T,1),[(-1)*ones(N-1,1) eye(N-1)]');
    Da=D*alpha;
    B=x1.*beta1+x2.*beta2+x3.*beta3+x4.*beta4;
    
    error=randn(N*T,1);
  %  error=-sqrt(3)+2*sqrt(3)*rand(N*T,1);
  %  error=(1/2)*chi2rnd(2,N*T,1)-1;
        
    Y=(eye(N*T)-rho_NT*W)\(B+Da+error);
    h=std(taut)*(T*N)^(-1/5);  %rule-of-thumb bandwidth

%%%% instrumental variables
H=[X W*X(:,2:end) W^2*X(:,2:end)];
Y_W=W*Y;

p=4;
S_H=zeros(N*T,N*T);
W_htau0=zeros(N*T,N*T);
Wstar_htau0=zeros(N*T,N*T);
for ii=1:T
    tau0=taut(ii,1);
for iii=1:N*T
    taoc0=tau(iii,1)-tau0;
    taoc(iii,:)=taoc0;% NT*1
end
    taoch=taoc./h; 
    M_Htau0=[H taoch.*H]; 
    W_htau0=diag(exp(-(taoch.^2)./2)/sqrt(2*pi)); 
    Wstar_htau0=(eye(N*T)-D*inv(D'*W_htau0*D)*D'*W_htau0)'*W_htau0*
    (eye(N*T)-D*inv(D'*W_htau0*D)*D'*W_htau0);
for iiii=1:N 
    sssH=[H(N*(ii-1)+iiii,:) zeros(1,3*p-2)]*inv(M_Htau0'*Wstar_htau0*M_Htau0)*M_Htau0'*Wstar_htau0;
    ssH(iiii,:)=sssH;  
end
    S_H(1+(ii-1)*N:ii*N,:)=ssH; 
end    
Bhat_H=S_H*Y_W;
psihat=inv(D'*D)*D'*(Y_W-Bhat_H);
Y_What=Bhat_H+D*psihat;

q=2;
Z_v=[Y_What X_v];
S_PL=zeros(N*T,N*T);
for i3=1:T
    tau0=taut(i3,1);
for i4=1:N*T
    taoc0=tau(i4,1)-tau0;
    taoc(i4,:)=taoc0;
end
    taoch=taoc./h; 
    M_PLtau0=[Z_v taoch.*Z_v];
    W_htau0=diag(exp(-(taoch.^2)./2)/sqrt(2*pi)); 
    Wstar_htau0=(eye(N*T)-D*inv(D'*W_htau0*D)*D'*W_htau0)'*W_htau0*
    (eye(N*T)-D*inv(D'*W_htau0*D)*D'*W_htau0);
    phi_PLtau0(1+(i3-1)*(q+1):i3*(q+1),:)=[eye(q+1) zeros(q+1,q+1)]*inv(M_PLtau0'*Wstar_htau0*M_PLtau0)*M_PLtau0'*Wstar_htau0;
for i5=1:N
    sssPL=[Z_v(N*(i3-1)+i5,:) zeros(1,q+1)]*inv(M_PLtau0'*Wstar_htau0*M_PLtau0)*M_PLtau0'*Wstar_htau0;
    ssPL(i5,:)=sssPL;  
end
    S_PL(1+(i3-1)*N:i3*N,:)=ssPL;  
end

Y_PLwave=(eye(N*T)-S_PL)*Y;
X_cwave=(eye(N*T)-S_PL)*X_c;
P_D=D*inv(D'*D)*D';
Y_PLbar=(eye(N*T)-P_D)*Y_PLwave;
X_cbar=(eye(N*T)-P_D)*X_cwave;
beta_chat=inv(X_cbar'*X_cbar)*X_cbar'*Y_PLbar;
beta3hat(i)=beta_chat(1,:);
beta4hat(i)=beta_chat(2,:);
gamma_vhat= phi_PLtau0*(Y-X_c*beta_chat);

rho00=[];
beta11=[];
beta22=[];
for j = 0:(length(gamma_vhat)/(q+1) - 1)
    rho00=[rho00,gamma_vhat((q+1)*j+ 1)];
    beta11=[beta11,gamma_vhat((q+1)*j+ 2)];
    beta22=[beta22,gamma_vhat((q+1)*j+ 3)];
end
rhohat(i,:)=rho00;
beta1hat(i,:)=beta11;
beta2hat(i,:)=beta22;
end  %for i=1:nsim  

%%%%%%%%%%%% Display results %%%%%%%%%%%%%%%%
T
N
rho_true=rhot';
rho_est=mean(rhohat);
rho_AMSE=mean(mean((rhohat-rho_true).^2))

beta1_true=beta1t'
beta1_est=mean(beta1hat)
beta1_AMSE=mean(mean((beta1hat-beta1_true).^2))

beta2_true=beta2t';
beta2_est=mean(beta2hat);
beta2_AMSE=mean(mean((beta2hat-beta2_true).^2))

beta3_est=mean(beta3hat);
beta3_bias=mean(beta3hat-beta3)
beta3_sd=std(beta3hat)

beta4_est=mean(beta4hat);
beta4_bias=mean(beta4hat-beta4)
beta4_sd=std(beta4hat)

figure;
plot(taut,rho_true,'-k','LineWidth',1);
hold on;
plot(taut,rho_est,'--r','LineWidth',1);
legend('True values', 'Average estimates'); 

figure;
plot(taut,beta1_true,'-k','LineWidth',1);
hold on;
plot(taut,beta1_est,'--r','LineWidth',1);
legend('True values', 'Average estimates'); 

figure;
plot(taut,beta2_true,'-k','LineWidth',1);
hold on;
plot(taut,beta2_est,'--r','LineWidth',1);
legend('True values', 'Average estimates'); 

disp("program ends");
diary off;
\end{lstlisting}

\subsection*{Testing}
\textbf{1. Size}
\begin{lstlisting}
function betat_size
%The test size for checking the linear relationship between some of explanatory variables and response variable

diary './rhot_size:T=3N=64error=1rook.txt';
disp("program starts");

%initialization
T=3;    %length of time
%T=5;
m=8;     %size of regular lattice
%m=10; 
N=m^2;    %number of regular points
nsim=100;  %number of simulations
nb=100;    %size of bootstrap sample


alpha0=[0.01 0.05 0.10];  %significance level
c=0;   %the degree of deviation from the null model

rand('seed',0);
randn('seed',1);

%generating row_standardized spatial adjacency matrix
w=generationrook(m,N);
%w=generationqueen(m,N);

sumc=sum(w,2); %sum by row
dd=ones(N,1)./sumc;
d=diag(dd);
ww=d*w;   %row normalization
W=kron(eye(T),ww);

number=zeros(1,3); %number of rejecting the null hypothesis
tt2=zeros(nsim,1);  %number of t*>t

%begin simulation
for i=1:nsim
        
    taut = (1:T)'./T;
    tau=kron(taut,ones(N,1));
    
    x1=ones(N*T,1);   %NT*1
    x2=randn(N*T,1);
    x3=randn(N*T,1);
    x4=randn(N*T,1);
    X=[x1 x2 x3 x4];  %NT*4
    X_v=[x1 x2];
    X_c=[x3 x4];
    
    rhot=-0.6*(sin(2*pi*taut).^2);
    %rhot=0.6*(sin(2*pi*taut).^2);
    rho_NT=kron(diag(rhot),eye(N));
    beta1t=4*taut;
    beta1=4*tau;
    beta2t=(taut+1).^2;
    beta2=(tau+1).^2;
    beta3t=-5+c*exp(taut);
    beta3=-5+c*exp(tau);
    beta4t=5+c*sin(2*pi*taut);  %T*1
    beta4=5+c*sin(2*pi*tau);  %nT*1
  
    alpha=rand(N-1,1); 
    D=kron(ones(T,1),[(-1)*ones(N-1,1) eye(N-1)]');
    Da=D*alpha;
    B=x1.*beta1+x2.*beta2+x3.*beta3+x4.*beta4;
    
    error=randn(N*T,1);
   % error=-sqrt(3)+2*sqrt(3)*rand(N*T,1);
    %error=(1/2)*chi2rnd(2,N*T,1)-1;
        
    Y=(eye(N*T)-rho_NT*W)\(B+Da+error);
    h=std(taut)*(T*N)^(-1/5);  %rule-of-thumb bandwidth

%%%% instrumental variables
H=[X W*X(:,2:end) W^2*X(:,2:end)];
Y_W=W*Y;

p=4;
q=2;
S_H=zeros(N*T,N*T);
W_htau0=zeros(N*T,N*T);
Wstar_htau0=zeros(N*T,N*T);
for ii=1:T
    tau0=taut(ii,1);
for iii=1:N*T
    taoc0=tau(iii,1)-tau0;
    taoc(iii,:)=taoc0;% NT*1
end
    taoch=taoc./h; 
    M_Htau0=[H taoch.*H]; 
    W_htau0=diag(exp(-(taoch.^2)./2)/sqrt(2*pi)); 
    Wstar_htau0=(eye(N*T)-D*inv(D'*W_htau0*D)*D'*W_htau0)'*W_htau0*
    (eye(N*T)-D*inv(D'*W_htau0*D)*D'*W_htau0);
    %phi_Htau0=[eye(3*p-2) zeros(3*p-2,3*p-2)]*inv(M_Htau0'*Wstar_htau0*M_Htau0)*M_Htau0'*Wstar_htau0;
for iiii=1:N 
    sssH=[H(N*(ii-1)+iiii,:) zeros(1,3*p-2)]*inv(M_Htau0'*Wstar_htau0*M_Htau0)*M_Htau0'*Wstar_htau0;
    ssH(iiii,:)=sssH;  
end
    S_H(1+(ii-1)*N:ii*N,:)=ssH; 
end    
Bhat_H=S_H*Y_W;
psihat=inv(D'*D)*D'*(Y_W-Bhat_H);
Y_What=Bhat_H+D*psihat;

%%%%%%%%%%% H0
Z_v=[Y_What X_v];
S_PL=zeros(N*T,N*T);
W_htau0=zeros(N*T,N*T);
Wstar_htau0=zeros(N*T,N*T);
for i3=1:T
    tau0=taut(i3,1);
for i4=1:N*T
    taoc0=tau(i4,1)-tau0;
    taoc(i4,:)=taoc0;
end
    taoch=taoc./h; 
    M_PLtau0=[Z_v taoch.*Z_v];
    W_htau0=diag(exp(-(taoch.^2)./2)/sqrt(2*pi)); 
    Wstar_htau0=(eye(N*T)-D*inv(D'*W_htau0*D)*D'*W_htau0)'*W_htau0*
    (eye(N*T)-D*inv(D'*W_htau0*D)*D'*W_htau0);
    phi_PLtau0(1+(i3-1)*(q+1):i3*(q+1),:)=[eye(q+1) zeros(q+1,q+1)]*inv(M_PLtau0'*Wstar_htau0*M_PLtau0)*M_PLtau0'*Wstar_htau0;
for i5=1:N
    sssPL=[Z_v(N*(i3-1)+i5,:) zeros(1,q+1)]*inv(M_PLtau0'*Wstar_htau0*M_PLtau0)*M_PLtau0'*Wstar_htau0;
    ssPL(i5,:)=sssPL;  
end
    S_PL(1+(i3-1)*N:i3*N,:)=ssPL;  
end

Y_PLwave=(eye(N*T)-S_PL)*Y;
X_cwave=(eye(N*T)-S_PL)*X_c;
P_D=D*inv(D'*D)*D';
Y_PLbar=(eye(N*T)-P_D)*Y_PLwave;
X_cbar=(eye(N*T)-P_D)*X_cwave;
beta_chat=inv(X_cbar'*X_cbar)*X_cbar'*Y_PLbar;
gamma_vhat= phi_PLtau0*(Y-X_c*beta_chat);

rhothat_pl=[];
beta1hat_pl=[];
beta2hat_pl=[];
for j = 0:(length(gamma_vhat)/(q+1) - 1)
    rhothat_pl=[rhothat_pl,gamma_vhat((q+1)*j+ 1)];
    beta1hat_pl=[beta1hat_pl,gamma_vhat((q+1)*j+ 2)];
    beta2hat_pl=[beta2hat_pl,gamma_vhat((q+1)*j+ 3)];
end
rhothat=rhothat_pl;
beta1hat=beta1hat_pl;
beta2hat=beta2hat_pl;
rhohat_NT=kron(diag(rhothat),eye(N));

%Bhat_Z_v=S_PL*(Y-X_c*beta_chat);
alpha_plhat=inv(D'*D)*D'*(Y_PLwave-X_cwave*beta_chat);
L_h=X_cbar*inv(X_cbar'*X_cbar)*X_cbar';
beta_rsd0=(eye(N*T)-L_h)*(eye(N*T)-P_D)*(eye(N*T)-S_PL)*Y;
beta_RSS0=beta_rsd0'*beta_rsd0;

%%%%%%%%%%%%%%%%% H1
Z=[Y_What X];

S_TV=zeros(N*T,N*T);
W_htau0=zeros(N*T,N*T);
Wstar_htau0=zeros(N*T,N*T);
for i6=1:T
    tau0=taut(i6,1);
for i7=1:N*T
    taoc0=tau(i7,1)-tau0;
    taoc(i7,:)=taoc0;
end
    taoch=taoc./h; 
    M_TVtau0=[Z taoch.*Z];
    W_htau0=diag(exp(-(taoch.^2)./2)/sqrt(2*pi)); 
    Wstar_htau0=(eye(N*T)-D*inv(D'*W_htau0*D)*D'*W_htau0)'*W_htau0*
    (eye(N*T)-D*inv(D'*W_htau0*D)*D'*W_htau0);
    %phi_TVtau0(1+(i6-1)*(p+1):i6*(p+1),:)=[eye(p+1) zeros(p+1,p+1)]*inv(M_TVtau0'*Wstar_htau0*M_TVtau0)*M_TVtau0'*Wstar_htau0;
for i8=1:N
    sssTV=[Z(N*(i6-1)+i8,:) zeros(1,p+1)]*inv(M_TVtau0'*Wstar_htau0*M_TVtau0)*M_TVtau0'*Wstar_htau0;
    ssTV(i8,:)=sssTV;  
end
    S_TV(1+(i6-1)*N:i6*N,:)=ssTV;  
end

beta_rsd1=(eye(N*T)-P_D)*(eye(N*T)-S_TV)*Y;
beta_RSS1=beta_rsd1'*beta_rsd1;

t=((N*T)/2)*((beta_RSS0-beta_RSS1)/beta_RSS1);

%bootstrapping the test statistic
epsilon=beta_rsd1;
epsibar=sum(epsilon)/(N*T);

%sampling bootstrap samples
for j=1:nb
estar=zeros(N*T,1);
   for jj=1:(N*T)
       rr=(N*T)*rand;
       intrr=round(rr);
       if intrr==0
          estar(jj,1)=epsilon(1,1)-epsibar;
       else
          estar(jj,1)=epsilon(intrr,1)-epsibar;
       end
   end

Ystar=(eye(N*T)-rhohat_NT*W)\(x1.*kron(beta1hat',ones(N,1))+x2.*
kron(beta2hat',ones(N,1))+X_c*beta_chat+D*alpha_plhat+estar);

%%%%%%instrumental variables
Y_Wstar=W*Ystar;
Bhatstar_H=S_H*Y_Wstar;
psihatstar=inv(D'*D)*D'*(Y_Wstar-Bhatstar_H);
Y_Whatstar=Bhatstar_H+D*psihatstar;

%H0
Z_vstar=[Y_Whatstar X_v];
for i3=1:T
    tau0=taut(i3,1);
for i4=1:N*T
    taoc0=tau(i4,1)-tau0;
    taoc(i4,:)=taoc0;
end
    taoch=taoc./h; 
    Mstar_PLtau0=[Z_vstar taoch.*Z_vstar];
    W_htau0=diag(exp(-(taoch.^2)./2)/sqrt(2*pi)); 
    Wstar_htau0=(eye(N*T)-D*inv(D'*W_htau0*D)*D'*W_htau0)'*W_htau0*
    (eye(N*T)-D*inv(D'*W_htau0*D)*D'*W_htau0);
for i5=1:N
    sssPLstar=[Z_vstar(N*(i3-1)+i5,:) zeros(1,q+1)]*inv(Mstar_PLtau0'*Wstar_htau0*Mstar_PLtau0)
    *Mstar_PLtau0'*Wstar_htau0;
    ssPLstar(i5,:)=sssPLstar;  
end
    S_PLstar(1+(i3-1)*N:i3*N,:)=ssPLstar;  
end

X_cwavestar=(eye(N*T)-S_PLstar)*X_c;
X_cbarstar=(eye(N*T)-P_D)*X_cwavestar;

L_hstar=X_cbarstar*inv(X_cbarstar'*X_cbarstar)*X_cbarstar';
beta_brsd0=(eye(N*T)-L_hstar)*(eye(N*T)-P_D)*(eye(N*T)-S_PLstar)*Ystar;
beta_bRSS0=beta_brsd0'*beta_brsd0;

%H1
Zstar=[Y_Whatstar X];

S_TVstar=zeros(N*T,N*T);
W_htau0=zeros(N*T,N*T);
Wstar_htau0=zeros(N*T,N*T);
for i9=1:T
    tau0=taut(i9,1);
for i10=1:N*T
    taoc0=tau(i10,1)-tau0;
    taoc(i10,:)=taoc0;
end
    taoch=taoc./h; 
    Mstar_TVtau0=[Zstar taoch.*Zstar];
    W_htau0=diag(exp(-(taoch.^2)./2)/sqrt(2*pi)); 
    Wstar_htau0=(eye(N*T)-D*inv(D'*W_htau0*D)*D'*W_htau0)'*W_htau0*
    (eye(N*T)-D*inv(D'*W_htau0*D)*D'*W_htau0);
for i11=1:N
    sssTVstar=[Zstar(N*(i9-1)+i11,:) zeros(1,p+1) ]*inv(Mstar_TVtau0'*Wstar_htau0*Mstar_TVtau0)*Mstar_TVtau0'*Wstar_htau0;
    ssTVstar(i11,:)=sssTVstar; 
end
    S_TVstar(1+(i9-1)*N:i9*N,:)=ssTVstar;  
end
beta_brsd1=(eye(N*T)-P_D)*(eye(N*T)-S_TVstar)*Ystar;
beta_bRSS1=beta_brsd1'*beta_brsd1;

tstar=((N*T)/2)*((beta_bRSS0-beta_bRSS1)/beta_bRSS1);

if tstar>t
   tt2(i)=tt2(i)+1;
else
   tt2(i)=tt2(i);
end
end %for j=1:nb


for i12=1:3
    if (tt2(i)/nb)<alpha0(1,i12)
      number(1,i12)=number(1,i12)+1;
    else
      number(1,i12)=number(1,i12);
    end
end

end  %for i=1:nsim 

%%%%%%%%%%% result
T
N
c
size=number/nsim

disp("program ends");
diary off;
\end{lstlisting}

\textbf{2. Power}
\begin{lstlisting}
function betat_power
%The test power for checking the linear relationship between some of explanatory variables and response variable

diary './rhot_size:T=3N=64error=1rook.txt';
disp("program starts");

%initialization
T=5;    %length of time
%T=5;
m=10;     %size of regular lattice
%m=10; 
N=m^2;    %number of regular points
nsim=500;  %number of simulations
nb=500;    %size of bootstrap sample


alpha0=0.05; %significance level
c0=[0.3 0.4 0.5];   %the degree of deviation from the null model

for c_i=1:3
c=c0(c_i);

rand('seed',0);
randn('seed',1);

%generating row_standardized spatial adjacency matrix
w=generationrook(m,N);
%w=generationqueen(m,N);

sumc=sum(w,2); %sum by row
dd=ones(N,1)./sumc;
d=diag(dd);
ww=d*w;   %row normalization
W=kron(eye(T),ww);

number=0; %number of rejecting the null hypothesis
tt2=zeros(nsim,1);  %number of t*>t

%begin simulation
for i=1:nsim
        
    taut = (1:T)'./T;
    tau=kron(taut,ones(N,1));
    
    x1=ones(N*T,1);   %NT*1
    x2=randn(N*T,1);
    x3=randn(N*T,1);
    x4=randn(N*T,1);
    X=[x1 x2 x3 x4];  %NT*4
    X_v=[x1 x2];
    X_c=[x3 x4];
    
    rhot=-0.6*(sin(2*pi*taut).^2);
    %rhot=0.6*(sin(2*pi*taut).^2);
    rho_NT=kron(diag(rhot),eye(N));
    beta1t=4*taut;
    beta1=4*tau;
    beta2t=(taut+1).^2;
    beta2=(tau+1).^2;
    beta3t=-5+c*exp(taut);
    beta3=-5+c*exp(tau);
    beta4t=5+c*sin(2*pi*taut);  %T*1
    beta4=5+c*sin(2*pi*tau);  %nT*1
  
    alpha=rand(N-1,1); 
    D=kron(ones(T,1),[(-1)*ones(N-1,1) eye(N-1)]');
    Da=D*alpha;
    B=x1.*beta1+x2.*beta2+x3.*beta3+x4.*beta4;
    
    %error=randn(N*T,1);
   error=-sqrt(3)+2*sqrt(3)*rand(N*T,1);
    %error=(1/2)*chi2rnd(2,N*T,1)-1;
        
    Y=(eye(N*T)-rho_NT*W)\(B+Da+error);
    h=std(taut)*(T*N)^(-1/5);  %rule-of-thumb bandwidth
    
%%%% instrumental variables
H=[X W*X(:,2:end) W^2*X(:,2:end)];
Y_W=W*Y;

p=4;
q=2;
S_H=zeros(N*T,N*T);
W_htau0=zeros(N*T,N*T);
Wstar_htau0=zeros(N*T,N*T);
for ii=1:T
    tau0=taut(ii,1);
for iii=1:N*T
    taoc0=tau(iii,1)-tau0;
    taoc(iii,:)=taoc0;% NT*1
end
    taoch=taoc./h; 
    M_Htau0=[H taoch.*H]; 
    W_htau0=diag(exp(-(taoch.^2)./2)/sqrt(2*pi)); 
    Wstar_htau0=(eye(N*T)-D*inv(D'*W_htau0*D)*D'*W_htau0)'*W_htau0*
    (eye(N*T)-D*inv(D'*W_htau0*D)*D'*W_htau0);
    %phi_Htau0=[eye(3*p-2) zeros(3*p-2,3*p-2)]*inv(M_Htau0'*Wstar_htau0*M_Htau0)*M_Htau0'*Wstar_htau0;
for iiii=1:N 
    sssH=[H(N*(ii-1)+iiii,:) zeros(1,3*p-2)]*inv(M_Htau0'*Wstar_htau0*M_Htau0)*M_Htau0'*Wstar_htau0;
    ssH(iiii,:)=sssH;  
end
    S_H(1+(ii-1)*N:ii*N,:)=ssH; 
end    
Bhat_H=S_H*Y_W;
psihat=inv(D'*D)*D'*(Y_W-Bhat_H);
Y_What=Bhat_H+D*psihat;

%%%%%%%%%%% H0
Z_v=[Y_What X_v];
S_PL=zeros(N*T,N*T);
W_htau0=zeros(N*T,N*T);
Wstar_htau0=zeros(N*T,N*T);
for i3=1:T
    tau0=taut(i3,1);
for i4=1:N*T
    taoc0=tau(i4,1)-tau0;
    taoc(i4,:)=taoc0;
end
    taoch=taoc./h; 
    M_PLtau0=[Z_v taoch.*Z_v];
    W_htau0=diag(exp(-(taoch.^2)./2)/sqrt(2*pi)); 
    Wstar_htau0=(eye(N*T)-D*inv(D'*W_htau0*D)*D'*W_htau0)'*W_htau0*
    (eye(N*T)-D*inv(D'*W_htau0*D)*D'*W_htau0);
    phi_PLtau0(1+(i3-1)*(q+1):i3*(q+1),:)=[eye(q+1) zeros(q+1,q+1)]*inv(M_PLtau0'*Wstar_htau0*M_PLtau0)*M_PLtau0'*Wstar_htau0;
for i5=1:N
    sssPL=[Z_v(N*(i3-1)+i5,:) zeros(1,q+1)]*inv(M_PLtau0'*Wstar_htau0*M_PLtau0)*M_PLtau0'*Wstar_htau0;
    ssPL(i5,:)=sssPL;  
end
    S_PL(1+(i3-1)*N:i3*N,:)=ssPL;  
end

Y_PLwave=(eye(N*T)-S_PL)*Y;
X_cwave=(eye(N*T)-S_PL)*X_c;
P_D=D*inv(D'*D)*D';
Y_PLbar=(eye(N*T)-P_D)*Y_PLwave;
X_cbar=(eye(N*T)-P_D)*X_cwave;
beta_chat=inv(X_cbar'*X_cbar)*X_cbar'*Y_PLbar;
gamma_vhat= phi_PLtau0*(Y-X_c*beta_chat);

rhothat_pl=[];
beta1hat_pl=[];
beta2hat_pl=[];
for j = 0:(length(gamma_vhat)/(q+1) - 1)
    rhothat_pl=[rhothat_pl,gamma_vhat((q+1)*j+ 1)];
    beta1hat_pl=[beta1hat_pl,gamma_vhat((q+1)*j+ 2)];
    beta2hat_pl=[beta2hat_pl,gamma_vhat((q+1)*j+ 3)];
end
rhothat=rhothat_pl;
beta1hat=beta1hat_pl;
beta2hat=beta2hat_pl;
rhohat_NT=kron(diag(rhothat),eye(N));

%Bhat_Z_v=S_PL*(Y-X_c*beta_chat);
alpha_plhat=inv(D'*D)*D'*(Y_PLwave-X_cwave*beta_chat);
L_h=X_cbar*inv(X_cbar'*X_cbar)*X_cbar';
beta_rsd0=(eye(N*T)-L_h)*(eye(N*T)-P_D)*(eye(N*T)-S_PL)*Y;
beta_RSS0=beta_rsd0'*beta_rsd0;

%%%%%%%%%%%%%%%%% H1
Z=[Y_What X];

S_TV=zeros(N*T,N*T);
W_htau0=zeros(N*T,N*T);
Wstar_htau0=zeros(N*T,N*T);
for i6=1:T
    tau0=taut(i6,1);
for i7=1:N*T
    taoc0=tau(i7,1)-tau0;
    taoc(i7,:)=taoc0;
end
    taoch=taoc./h; 
    M_TVtau0=[Z taoch.*Z];
    W_htau0=diag(exp(-(taoch.^2)./2)/sqrt(2*pi)); 
    Wstar_htau0=(eye(N*T)-D*inv(D'*W_htau0*D)*D'*W_htau0)'*W_htau0*
    (eye(N*T)-D*inv(D'*W_htau0*D)*D'*W_htau0);
    %phi_TVtau0(1+(i6-1)*(p+1):i6*(p+1),:)=[eye(p+1) zeros(p+1,p+1)]*inv(M_TVtau0'*Wstar_htau0*M_TVtau0)*M_TVtau0'*Wstar_htau0;
for i8=1:N
    sssTV=[Z(N*(i6-1)+i8,:) zeros(1,p+1)]*inv(M_TVtau0'*Wstar_htau0*M_TVtau0)*M_TVtau0'*Wstar_htau0;
    ssTV(i8,:)=sssTV;  
end
    S_TV(1+(i6-1)*N:i6*N,:)=ssTV;  
end

beta_rsd1=(eye(N*T)-P_D)*(eye(N*T)-S_TV)*Y;
beta_RSS1=beta_rsd1'*beta_rsd1;

t=((N*T)/2)*((beta_RSS0-beta_RSS1)/beta_RSS1);

%bootstrapping the test statistic
epsilon=beta_rsd1;
epsibar=sum(epsilon)/(N*T);

%sampling bootstrap samples
for j=1:nb
estar=zeros(N*T,1);
   for jj=1:(N*T)
       rr=(N*T)*rand;
       intrr=round(rr);
       if intrr==0
          estar(jj,1)=epsilon(1,1)-epsibar;
       else
          estar(jj,1)=epsilon(intrr,1)-epsibar;
       end
   end

Ystar=(eye(N*T)-rhohat_NT*W)\(x1.*kron(beta1hat',ones(N,1))+x2.*
kron(beta2hat',ones(N,1))+X_c*beta_chat+D*alpha_plhat+estar);

%%%%%%instrumental variables
Y_Wstar=W*Ystar;
Bhatstar_H=S_H*Y_Wstar;
psihatstar=inv(D'*D)*D'*(Y_Wstar-Bhatstar_H);
Y_Whatstar=Bhatstar_H+D*psihatstar;

%H0
Z_vstar=[Y_Whatstar X_v];
for i3=1:T
    tau0=taut(i3,1);
for i4=1:N*T
    taoc0=tau(i4,1)-tau0;
    taoc(i4,:)=taoc0;
end
    taoch=taoc./h; 
    Mstar_PLtau0=[Z_vstar taoch.*Z_vstar];
    W_htau0=diag(exp(-(taoch.^2)./2)/sqrt(2*pi)); 
    Wstar_htau0=(eye(N*T)-D*inv(D'*W_htau0*D)*D'*W_htau0)'*W_htau0*
    (eye(N*T)-D*inv(D'*W_htau0*D)*D'*W_htau0);
for i5=1:N
    sssPLstar=[Z_vstar(N*(i3-1)+i5,:) zeros(1,q+1)]*inv(Mstar_PLtau0'*Wstar_htau0*Mstar_PLtau0)*Mstar_PLtau0'*
    Wstar_htau0;
    ssPLstar(i5,:)=sssPLstar;  
end
    S_PLstar(1+(i3-1)*N:i3*N,:)=ssPLstar;  
end

X_cwavestar=(eye(N*T)-S_PLstar)*X_c;
X_cbarstar=(eye(N*T)-P_D)*X_cwavestar;

L_hstar=X_cbarstar*inv(X_cbarstar'*X_cbarstar)*X_cbarstar';
beta_brsd0=(eye(N*T)-L_hstar)*(eye(N*T)-P_D)*(eye(N*T)-S_PLstar)*Ystar;
beta_bRSS0=beta_brsd0'*beta_brsd0;

%H1
Zstar=[Y_Whatstar X];

S_TVstar=zeros(N*T,N*T);
W_htau0=zeros(N*T,N*T);
Wstar_htau0=zeros(N*T,N*T);
for i9=1:T
    tau0=taut(i9,1);
for i10=1:N*T
    taoc0=tau(i10,1)-tau0;
    taoc(i10,:)=taoc0;
end
    taoch=taoc./h; 
    Mstar_TVtau0=[Zstar taoch.*Zstar];
    W_htau0=diag(exp(-(taoch.^2)./2)/sqrt(2*pi)); 
    Wstar_htau0=(eye(N*T)-D*inv(D'*W_htau0*D)*D'*W_htau0)'*W_htau0*
    (eye(N*T)-D*inv(D'*W_htau0*D)*D'*W_htau0);
for i11=1:N
    sssTVstar=[Zstar(N*(i9-1)+i11,:) zeros(1,p+1) ]*inv(Mstar_TVtau0'*Wstar_htau0*Mstar_TVtau0)*Mstar_TVtau0'*Wstar_htau0;
    ssTVstar(i11,:)=sssTVstar; 
end
    S_TVstar(1+(i9-1)*N:i9*N,:)=ssTVstar;  
end
beta_brsd1=(eye(N*T)-P_D)*(eye(N*T)-S_TVstar)*Ystar;
beta_bRSS1=beta_brsd1'*beta_brsd1;

tstar=((N*T)/2)*((beta_bRSS0-beta_bRSS1)/beta_bRSS1);

if tstar>t
   tt2(i)=tt2(i)+1;
else
   tt2(i)=tt2(i);
end
end %for j=1:nb


if (tt2(i)/nb)<alpha0
    number=number+1;
else
    number=number;
end

end  %for i=1:nsim 

%%%%%%%%%%% result
T
N
c
power=number/nsim
end

disp("program ends");
diary off;
\end{lstlisting}

\subsection*{Application}
\begin{lstlisting}
function betat_test_data
%The test for time-varying coefficients in Carbon emission data

diary './rho_test_data.txt';
disp("program starts");

%initialization
T=12;    %length of time
N=30;    %number of regular points
nb=500;    %size of bootstrap sample
tt2=0;  %number of t*>t

%loading data
data=xlsread('data.xlsx');
PC=data(:,1); 
PG=data(:,2);  
PR=data(:,3);   
IR=data(:,4);    
ER=data(:,5);

% read weight matrix
w=readtable('wmatrix.csv');
w_matrix = table2array(w);
sumc=sum(w_matrix,2); %sum by row
dd=ones(N,1)./sumc;
d=diag(dd);
ww=d*w_matrix;   %row normalization
W=kron(eye(T),ww);

taut = (1:T)'./T;
tau=kron(taut,ones(N,1));
h=std(taut)*(T*N)^(-1/5);  %rule-of-thumb bandwidth
D=kron(ones(T,1),[(-1)*ones(N-1,1) eye(N-1)]');

X=[ones(N*T,1) PG PR IR ER];
Y=PC;
 
X_v=[ones(N*T,1) PR ER];
X_c=[PG IR];

%%%% instrumental variables
H=[X W*X(:,2:end) W^2*X(:,2:end)];
Y_W=W*Y;

p=5;
S_H=zeros(N*T,N*T);
for ii=1:T
    tau0=taut(ii,1);
for iii=1:N*T
    taoc0=tau(iii,1)-tau0;
    taoc(iii,:)=taoc0;% NT*1
end
    taoch=taoc./h; 
    M_Htau0=[H taoch.*H]; 
    W_htau0=diag(exp(-(taoch.^2)./2)/sqrt(2*pi)); 
    Wstar_htau0=(eye(N*T)-D*inv(D'*W_htau0*D)*D'*W_htau0)'*W_htau0*
    (eye(N*T)-D*inv(D'*W_htau0*D)*D'*W_htau0);
for iiii=1:N 
    sssH=[H(N*(ii-1)+iiii,:) zeros(1,3*p-2)]*inv(M_Htau0'*Wstar_htau0*M_Htau0)*M_Htau0'*Wstar_htau0;
    ssH(iiii,:)=sssH;  
end
    S_H(1+(ii-1)*N:ii*N,:)=ssH; 
end    
Bhat_H=S_H*Y_W;
psihat=inv(D'*D)*D'*(Y_W-Bhat_H);
Y_What=Bhat_H+D*psihat;

%%%%%%%%%%% H0
q=3;
Z_v=[Y_What X_v];
S_PL=zeros(N*T,N*T);
for i3=1:T
    tau0=taut(i3,1);
for i4=1:N*T
    taoc0=tau(i4,1)-tau0;
    taoc(i4,:)=taoc0;
end
    taoch=taoc./h; 
    M_PLtau0=[Z_v taoch.*Z_v];
    W_htau0=diag(exp(-(taoch.^2)./2)/sqrt(2*pi)); 
    Wstar_htau0=(eye(N*T)-D*inv(D'*W_htau0*D)*D'*W_htau0)'*W_htau0*
    (eye(N*T)-D*inv(D'*W_htau0*D)*D'*W_htau0);
    phi_PLtau0(1+(i3-1)*(q+1):i3*(q+1),:)=[eye(q+1) zeros(q+1,q+1)]*inv(M_PLtau0'*Wstar_htau0*M_PLtau0)*M_PLtau0'*Wstar_htau0;
for i5=1:N
    sssPL=[Z_v(N*(i3-1)+i5,:) zeros(1,q+1)]*inv(M_PLtau0'*Wstar_htau0*M_PLtau0)*M_PLtau0'*Wstar_htau0;
    ssPL(i5,:)=sssPL;  
end
    S_PL(1+(i3-1)*N:i3*N,:)=ssPL;  
end

Y_PLwave=(eye(N*T)-S_PL)*Y;
X_cwave=(eye(N*T)-S_PL)*X_c;
P_D=D*inv(D'*D)*D';
Y_PLbar=(eye(N*T)-P_D)*Y_PLwave;
X_cbar=(eye(N*T)-P_D)*X_cwave;

beta_chat=inv(X_cbar'*X_cbar)*X_cbar'*Y_PLbar;
beta3hat=beta_chat(1,:);
beta4hat=beta_chat(2,:);

gamma_vhat= phi_PLtau0*(Y-X_c*beta_chat);
rhohat_pl=[];
beta0hat_pl=[];
beta1hat_pl=[];
beta2hat_pl=[];
for j = 0:(length(gamma_vhat)/(q+1) - 1)
    rhohat_pl=[rhohat_pl,gamma_vhat((q+1)*j+ 1)];
    beta0hat_pl=[beta0hat_pl,gamma_vhat((q+1)*j+ 2)];
    beta1hat_pl=[beta1hat_pl,gamma_vhat((q+1)*j+ 3)];
    beta2hat_pl=[beta2hat_pl,gamma_vhat((q+1)*j+ 4)];
end
rhothat=rhohat_pl;
beta0hat=beta0hat_pl;
beta1hat=beta1hat_pl;
beta2hat=beta2hat_pl;

%Bhat_Z_v=S_PL*(Y-X_c*beta_chat);
alpha_plhat=inv(D'*D)*D'*(Y_PLwave-X_cwave*beta_chat);
L_PL=X_cbar*inv(X_cbar'*X_cbar)*X_cbar';
beta_rsd0=(eye(N*T)-L_PL)*(eye(N*T)-P_D)*(eye(N*T)-S_PL)*Y;
beta_RSS0=beta_rsd0'*beta_rsd0;

%%%%%%%%%%%%%%%%% H1
Z=[Y_What X];
S_TV=zeros(N*T,N*T);
for i6=1:T
    tau0=taut(i6,1);
for i7=1:N*T
    taoc0=tau(i7,1)-tau0;
    taoc(i7,:)=taoc0;
end
    taoch=taoc./h; 
    M_TVtau0=[Z taoch.*Z];
    W_htau0=diag(exp(-(taoch.^2)./2)/sqrt(2*pi)); 
    Wstar_htau0=(eye(N*T)-D*inv(D'*W_htau0*D)*D'*W_htau0)'*W_htau0*
    (eye(N*T)-D*inv(D'*W_htau0*D)*D'*W_htau0);
    %phi_TVtau0(1+(i6-1)*(p+1):i6*(p+1),:)=[eye(p+1) zeros(p+1,p+1)]*inv(M_TVtau0'*Wstar_htau0*M_TVtau0)*M_TVtau0'*Wstar_htau0;
    %gammahat(1+(i6-1)*(p+1):i6*(p+1),:)=phi_TVtau0(1+(i6-1)*(p+1):i6*(p+1),:)*Y;
for i8=1:N
    sssTV=[Z(N*(i6-1)+i8,:) zeros(1,p+1) ]*inv(M_TVtau0'*Wstar_htau0*M_TVtau0)*M_TVtau0'*Wstar_htau0;
    ssTV(i8,:)=sssTV;  
end
    S_TV(1+(i6-1)*N:i6*N,:)=ssTV;  
end

beta_rsd1=(eye(N*T)-P_D)*(eye(N*T)-S_TV)*Y;
beta_RSS1=beta_rsd1'*beta_rsd1;

t=((N*T)/2)*((beta_RSS0-beta_RSS1)/beta_RSS1);

%bootstrapping the test statistic
epsilon=beta_rsd1;
epsibar=sum(epsilon)/(N*T);

rand('seed',555);
for j=1:nb
%sampling bootstrap samples
estar=zeros(N*T,1);
   for jj=1:(N*T)
       rr=(N*T)*rand;
       intrr=round(rr);
       if intrr==0
          estar(jj,1)=epsilon(1,1)-epsibar;
       else
          estar(jj,1)=epsilon(intrr,1)-epsibar;
       end
   end 

rhohat_NT=kron(diag(rhothat),eye(N));
Ystar=(eye(N*T)-rhohat_NT*W)\(ones(N*T,1).*kron(beta0hat',ones(N,1))+PR.*kron
(beta1hat',ones(N,1))+ER.*kron(beta2hat',ones(N,1))+X_c*beta_chat+D*alpha_plhat+estar);

%%%%%%instrumental variables
Y_Wstar=W*Ystar;
Bhatstar_H=S_H*Y_Wstar;
psihatstar=inv(D'*D)*D'*(Y_Wstar-Bhatstar_H);
Y_Whatstar=Bhatstar_H+D*psihatstar;

%H0
Z_vstar=[Y_Whatstar X_v];
S_PLstar=zeros(N*T,N*T);
for i3=1:T
    tau0=taut(i3,1);
for i4=1:N*T
    taoc0=tau(i4,1)-tau0;
    taoc(i4,:)=taoc0;
end
    taoch=taoc./h; 
    Mstar_PLtau0=[Z_vstar taoch.*Z_vstar];
    W_htau0=diag(exp(-(taoch.^2)./2)/sqrt(2*pi)); 
    Wstar_htau0=(eye(N*T)-D*inv(D'*W_htau0*D)*D'*W_htau0)'*W_htau0*
    (eye(N*T)-D*inv(D'*W_htau0*D)*D'*W_htau0);
for i5=1:N
    sssPLstar=[Z_vstar(N*(i3-1)+i5,:) zeros(1,q+1)]*inv(Mstar_PLtau0'*Wstar_htau0*Mstar_PLtau0)*
    Mstar_PLtau0'*Wstar_htau0;
    ssPLstar(i5,:)=sssPLstar;  
end
    S_PLstar(1+(i3-1)*N:i3*N,:)=ssPLstar;  
end
X_cwavestar=(eye(N*T)-S_PLstar)*X_c;
X_cbarstar=(eye(N*T)-P_D)*X_cwavestar;
L_PLstar=X_cbarstar*inv(X_cbarstar'*X_cbarstar)*X_cbarstar';
beta_brsd0=(eye(N*T)-L_PLstar)*(eye(N*T)-P_D)*(eye(N*T)-S_PLstar)*Ystar;
beta_bRSS0=beta_brsd0'*beta_brsd0;

%H1
Zstar=[Y_Whatstar X];
S_TVstar=zeros(N*T,N*T);
for i9=1:T
    tau0=taut(i9,1);
for i10=1:N*T
    taoc0=tau(i10,1)-tau0;
    taoc(i10,:)=taoc0;
end
    taoch=taoc./h; 
    Mstar_TVtau0=[Zstar taoch.*Zstar];
    W_htau0=diag(exp(-(taoch.^2)./2)/sqrt(2*pi)); 
    Wstar_htau0=(eye(N*T)-D*inv(D'*W_htau0*D)*D'*W_htau0)'*W_htau0*
    (eye(N*T)-D*inv(D'*W_htau0*D)*D'*W_htau0);
for i11=1:N
    sssTVstar=[Zstar(N*(i9-1)+i11,:) zeros(1,p+1) ]*inv(Mstar_TVtau0'*Wstar_htau0*Mstar_TVtau0)*Mstar_TVtau0'*Wstar_htau0;
    ssTVstar(i11,:)=sssTVstar; 
end
    S_TVstar(1+(i9-1)*N:i9*N,:)=ssTVstar;  
end
rho_brsd1=(eye(N*T)-P_D)*(eye(N*T)-S_TVstar)*Ystar;
rho_bRSS1=rho_brsd1'*rho_brsd1;

tstar=((N*T)/2)*((beta_bRSS0-rho_bRSS1)/rho_bRSS1);

if tstar>t
   tt2=tt2+1;
else
   tt2=tt2;
end
end %for j=1:nb

%%%%%%%%%%% result
p=tt2/nb

beta3hat
beta4hat

rhothat;
beta0hat;
beta1hat;
beta2hat;

% Plotting rhothat, beta0hat, beta1hat, beta2hat
year=2005:2016;
figure;
plot(year,rhothat,'k','LineWidth',1);
xlabel('Time');
ylabel('rhothat');
figure;
plot(year,beta0hat,'k','LineWidth',1);
xlabel('Time');
ylabel('beta0hat');
figure;
plot(year,beta1hat,'k','LineWidth',1);
xlabel('Time');
ylabel('beta1hat');
figure;
plot(year,beta2hat,'k','LineWidth',1);
xlabel('beta2hat');

disp("program ends");
diary off;
end
\end{lstlisting}

\end{document}